\documentclass[12pt]{article}
\usepackage{amsmath}
\newcommand{\be}{\begin{equation}}
\newcommand{\ee}{\end{equation}}
\newcommand{\ba}{\begin{eqnarray}}
\newcommand{\ea}{\end{eqnarray}}
\newcommand{\ban}{\begin{eqnarray*}}
\newcommand{\ean}{\end{eqnarray*}}

 \newcommand{\qed}{\hspace*{\fill}\rule{3mm}{3mm}\quad}
\newcommand{\Pf}{\noindent  {\em Proof.} }

\newcommand{\sect}[1]{\section{#1}  \setcounter{equation}{0}}

\newtheorem{lem}{Lemma}[section]
 
\begin{document}
\newtheorem{defn}[lem]{Definition}
\newtheorem{theo}[lem]{Theorem}
\newtheorem{cor}[lem]{Corollary}
\newtheorem{prop}[lem]{Proposition}
\newtheorem{rk}[lem]{Remark}
\newtheorem{ex}[lem]{Example}
\newtheorem{note}[lem]{Note}
\newtheorem{conj}[lem]{Conjecture}

\title{The Logarithmic Sobolev Inequality Along The Ricci Flow \\ 
(revised version) }
\author{Rugang Ye \\ {\small Department  of Mathematics} \\
{\small University of California, Santa  Barbara}}
\date{July 20, 2007}
\maketitle

\noindent 1. Introduction \\
2. The Sobolev inequality \\
3. The logarithmic Sobolev inequality on a Riemannian manifold \\
4. The logarithmic Sobolev inequality along the Ricci flow \\
5. The Sobolev inequality along the Ricci flow\\
6. The $\kappa$-noncollapsing estimate \\
Appendix A. The logarithmic Sobolev inequalities on the euclidean space \\
Appendix B. The estimate of $e^{-tH}$ \\
Appendix C. From the estimate for $e^{-tH}$ to the Sobolev inequality \\

\sect{Introduction}
 
Consider a compact manifold $M$ of dimension $n \ge 3$.  Let $g=g(t)$ be a smooth solution of the Ricci flow 
\ba
\frac{\partial g}{\partial t}=-2Ric
\ea
on $M \times [0, T)$ for some (finite or infinite) $T>0$ with a given initial metric 
$g(0)=g_0$.\\

\noindent {\bf Theorem A} {\it For each $\sigma>0$ and each $t \in [0, T)$  there holds
\ba \label{sobolevA}
\int_M u^2 \ln u^2 dvol &\le& \sigma \int_M (|\nabla u|^2 +\frac{R}{4} u^2)dvol  -\frac{n}{2} \ln \sigma 
+A_1(t+\frac{\sigma}{4})+A_2
\ea
for all $u\in W^{1,2}(M)$ with  $\int_M u^2 dvol=1$, 
where
\ba
A_1&=&\frac{4}{\tilde C_S(M, g_0)^2 vol_{g_0}(M)^{\frac{2}{n}}}-\min R_{g_0}, \nonumber \\
A_2&=& n\ln \tilde C_S(M,g_0)+
\frac{n}{2}(\ln n-1), \nonumber
\ea
and 
all geometric quantities  are associated with the metric $g(t)$ (e.g. the volume form $dvol$ and the scalar curvature $R$), except the scalar curvature $R_{g_0}$, the modified Sobolev constant $\tilde C_S(M,g_0)$ (see Section 2 for its definition) and the volume $vol_{g_0}(M)$ which 
are those of the initial metric $g_0$. 

Consequently, there holds for each $t \in [0, T)$
\ba \label{strongA}
\int_M u^2 \ln u^2 dvol \le \frac{n}{2} \ln \left[ \alpha_I(\int_M (|\nabla u|^2 +\frac{R}{4} u^2)dvol+\frac{A_1}{4}) \right]
\ea
for all $u \in W^{1,2}(M)$ with $\int_M u^2 dvol=1$, 
where 
\ba
\alpha_I=\frac{2e}{n}e^{\frac{2(A_1t+A_2)}{n}}.
\ea
} \\

The exact factor $\frac{n}{2}$ in the term $-\frac{n}{2} \ln \sigma$ in the logarithmic Sobolev inequality (\ref{sobolevA}) (also in (\ref{sobolevB}) and (\ref{sobolevC}) below) is 
crucial for the purpose of Theorem D and Theorem $\mbox{D}^*$. Note that an upper bound for the Sobolev constant $C_S(M,g_0)$ and 
the modified Sobolev constant $\tilde C_S(M, g_0)$ can be obtained in terms of 
a lower bound for the diameter rescaled Ricci curvature and a positive lower bound for 
the diameter rescaled volume, see Section 2. In particular, a lower bound for the Ricci curvature,
a positive lower bound for the volume and an upper bound for the diameter lead to an upper bound for 
the Sobolev constant and the modified Sobolev constant.

The logarithmic Sobolev inequality in Theorem A is uniform for all time which lies below a given bound, but deteriorates as 
time becomes large. The next  result takes care of large time under the assumption that a certain  eigenvalue 
$\lambda_0$ of the initial metric is positive. This assumption holds true e.g. when the scalar curvature is 
nonnegative and somewhere positive. \\

\noindent {\bf Theorem B} {\it Assume that the first eigenvalue $\lambda_0=\lambda_0(g_0)$ of the operator 
$-\Delta+\frac{R}{4}$ for  the initial metric $g_0$ is positive. Let $\delta_0=\delta_0(g_0)$ be the number defined in 
(\ref{delta-0}). Let $t \in [0, T)$ and $\sigma>0$ satisfy $t+\sigma \ge \frac{n}{8}C_S(M,g_0)^2\delta_0$.
Then there holds
\ba \label{sobolevB}
\int_M u^2 \ln u^2 dvol &\le& \sigma \int_M (|\nabla u|^2+\frac{R}{4}u^2)dvol -\frac{n}{2}\ln \sigma 
\nonumber \\
&&+\frac{n}{2}\ln n + n \ln C_S(M,g_0)+\sigma_0(g_0)
\ea
for all $u\in W^{1,2}(M)$ with  $\int_M u^2 dvol=1$, 
where all geometric quantities  are associated with the metric $g(t)$ (e.g. the volume form $dvol$ and the scalar curvature $R$), except the Sobolev constant $C_S(M,g_0)$ and the number 
$\sigma_0(g_0)$ (defined in (\ref{sigma-0}))  which 
are those of the initial metric $g_0$.

Consequently, there holds for each $t \in [0, T)$
\ba \label{strongB}
\int_M u^2 \ln u^2 dvol \le \frac{n}{2} \ln \left[ \alpha_{II}\int_M (|\nabla u|^2 +\frac{R}{4} u^2)dvol \right]
\ea
for all $u \in W^{1,2}(M)$ with $\int_M u^2 dvol=1$, 
where 
\ba
\alpha_{II}=2eC_S(M,g_0)^2 e^{\frac{2}{n}\sigma_0(g_0)}.
\ea
} \\

Combining Theorem A and Theorem B we obtain a uniform  logarithmic Sobolev inequality along the Ricci 
flow without any restriction on time or the factor $\sigma$, assuming only $\lambda_0(g_0) > 0$.  \\

\noindent {\bf Theorem C} {\it Assume that $\lambda_0(g_0) > 0$. For each $t\in [0, T)$ and each $\sigma>0$ there holds 
\ba \label{sobolevC}
\int_M u^2 \ln u^2 dvol \le \sigma \int_M (|\nabla u|^2 +\frac{R}{4}u^2) dvol
-\frac{n}{2}\ln \sigma +C
\ea
for all $u \in W^{1,2}(M)$ with $\int_M u^2 dvol=1$, where $C$ depends only on the dimension $n$, a positive lower bound for $vol_{g_0}(M)$, a nonpositive 
lower bound for $R_{g_0}$, an upper bound for $C_S(M, g_0)$, and a positive lower bound for 
$\lambda_0(g_0)$.

Consequently, there holds for each $t \in [0, T)$
\ba \label{strongC}
\int_M u^2 \ln u^2 dvol \le \frac{n}{2} \ln \left[ \alpha_{III}\int_M (|\nabla u|^2 +\frac{R}{4} u^2)dvol \right]
\ea
for all $u \in W^{1,2}(M)$ with $\int_M u^2 dvol=1$, 
where 
\ba
\alpha_{III}=\frac{2e}{n}e^{\frac{2}{n}C}.
\ea
} \\

We note here a special consequence of Theorem C. \\

\noindent {\bf Corollary} {\it  Assume that $\lambda_0(g_0)>  0$. Then we have at any time $t \in [0, T)$
\ba \label{volume3}
vol_{g(t)}(M) \ge e^{-\frac{1}{4}-C}
\ea
when $\hat R(t) \le 0$, and 
\ba \label{volume4}
vol_{g(t)}(M) \ge e^{-\frac{1}{4}-C} \hat R(t)^{-\frac{n}{2}}
\ea
when $\hat R(t)>0$. Here $\hat R$ denotes the average scalar curvature.} \\

Similar volume bounds follow from Theorem A without the condition $\lambda_0(g_0)> 0$, but they 
also depend on a (finite) upper bound of $T$.

The class of Riemannian manifolds $(M, g_0)$ with $\lambda_0(g_0) > 0$  (or, more generally, 
$\lambda_0(g_0)\ge 0$) is a very large one
and particularly significant from a geometric point of view. 
On the other hand, we  would like to point out that the assumption $\lambda_0(g_0)>0$ (or, more generally, 
$\lambda_0(g_0)\ge 0$) in 
Theorem C is indispensible in general. (The case 
$\lambda_0(g_0)=0$ will be presented elsewhere.)    In other words, a uniform logarithmic 
Sobolev inequality like (\ref{sobolevC}) without the assumption $\lambda_0(g_0) \ge 0$ 
is false in general. Indeed, by [HI] there are  smooth solutions of the Ricci flow 
on torus bundles over the circle which exist for all time, have  bounded curvature, 
and collapse as $t\rightarrow \infty$. In view of the proofs of Theorem D and Theorem E, a  uniform logarithmic Sobolev inequality like (\ref{sobolevC}) fails to hold along these solutions. The generalization of Theorem C stated in the first posted  version of 
this paper is thus incorrect. The trouble with the proof,  
which we found in the process of 
trying to work out an improvement of the logarithmic Sobolev inequality along the Ricci flow, stems from 
the application of a monotonicity formula in [Z].
More precisely, for a given $T^*>0$,  the estimate 
\ba \label{F}
F\equiv \int_M (R+|\nabla f|^2) \frac{e^{-f}}{(4\pi \tau)^{\frac{n}{2}}} dvol \le \frac{2n}{\tau}
\ea
 in [Z] (needed for proving the monotonicity of the 
generalized entropy in [Z]) holds true only for 
a solution $f$ of the equation (\ref{nonlinearconjugate}) defined up to $T^*$, where $\tau=T^*-t$. For the choice 
$T^*=2 t_0$ for a given $t_0$ in [Z]  and the first posted version of this paper, where $t_0$ is denoted 
$t$, $f$ is assumed to start backwards at $t_0$, and may not 
exist on $[t_0, T^*)$. Hence the inequality (\ref{F}) may not hold, and  
the generalized entropy may not be monotone.  In contrast, Perelman's entropy monotonicity is 
always valid in any time interval  where $f$ satisfies (\ref{nonlinearconjugate}). In other words, for a given $t_0$, one can choose 
$T^*> t_0$ arbitrarily to define $\tau=T^*-t$.  Perelman's entropy is monotone on $[0, t_0]$ as long as 
$f$ satisfies (\ref{nonlinearconjugate}) there.  This is crucial for applying
Perelman's entropy.
 
For a brief account of the logarithmic Sobolev inequalitities on the euclidean space we refer 
to Appendix A, which serve as the background for the idea of the logarithmic Sobolev inequality.  Both Theorem A and Theorem B are consequences of Perelman's entropy monotonicity [P1]. We obtained these  two results, Theorem C 
and Theorem 4.2 in 2004 (around the time of the author's differential geometry seminar talk ``An introduction 
to the logarithmic Sobolev inequality" at UCSB in June 2004).  They have also 
been prepared as part of the notes [Y4].  

Inspired by an argument in [Z], we apply the theory as presented in Chapter 2 of [D] 
to derive from Theorem D a Sobolev inequality along the Ricci flow without any 
restriction on time. A particularly nice feature of the theory in Chapter 2 of [D] is that no additional geometric data (such as the volume)
are involved in the passage from the logarithmic Sobolev inequality to the Sobolev inequality. Only 
the non-integral terms in the logarithmic Sobolev inequality and a nonpositive lower bound 
for the potential function $\Psi$ (see Theorem \ref{D3}) come into play. This leads to  
the form of the geometric dependenc in the following theorem. \\

\noindent {\bf Theorem D} {\it Assume that $\lambda_0(g_0) > 0$.  There is a  positive constant $A$ depending only on 
the dimension $n$, a nonpositive lower bound for $R_{g_0}$, a positive lower bound for $vol_{g_0}(M)$, an upper bound for $C_S(M,g_0)$,
and a positive lower bound for $\lambda_0(g_0)$, such that for each $t \in [0, T)$ and 
all $u \in W^{1,2}(M)$ there holds 
\ba \label{sobolevD}
\left( \int_M |u|^{\frac{2n}{n-2}} dvol \right)^{\frac{n-2}{n}} \le 
A\int_M (|\nabla u|^2+\frac{R}{4}u^2) dvol,
\ea
where all geometric quantities except $A$  are associated with $g(t)$. 
}\\

In a similar fashion, a Sobolev inequality follows from Theorem A in which the condition 
$\lambda_0(g_0) > 0$ is not assumed, but the bounds also depend on an upper bound of time. 
\\

\noindent {\bf Theorem $\mbox{D}^*$} {\it  Assume $T<\infty$. There are positive constants $A$ and $B$ depending only on the dimension $n$, 
a nonpositive lower bound for $R_{g_0}$, a positive lower bound for $vol_{g_0}(M)$, an upper bound for $C_S(M,g_0)$,
and an upper bound for $T$, such that for each $t \in [0, T)$ and 
all $u \in W^{1,2}(M)$ there holds 
\ba \label{sobolevD*}
\left( \int_M |u|^{\frac{2n}{n-2}} dvol \right)^{\frac{n-2}{n}} \le 
A\int_M (|\nabla u|^2+\frac{R}{4}u^2)dvol+B \int_M u^2 dvol,
\ea
where all geometric quantities except $A$ and $B$  are associated with $g(t)$. 
}\\

We also obtain two results which extend Theorem D and Theorem $\mbox{D}^*$ to 
the set-up of $W^{1,p}(M)$ for all $1<p<n$, see Theorem C.6 and Theorem C.7 in Appendix C. 
(Theorem D and Theorem $\mbox{D}^*$ correspond to the 
case $p=2$.)  These two general results can be thought of  as 
nonlocal versions of Sobolev inequality, because they involve nonlocal pseudo-differential operators.  Further results on (conventional) Sobolev inequalities
for $2<p<n$ and $1<p<2$ will be presented elsewhere, which are derived from Theorem D, Theorem $\mbox{D}^*$, Theorem C.6, Theorem C.7 and a result on Riesz transforms.  
 We would like to point out that 
the $p=2$ case of the Sobolev inequality is the most important for analytic and geometric applications.

The theory in Chapter 2 of [D] is formulated in a general and abstract set-up of symmetric Markov semigroups. 
By Lemma \ref{contractionlemma}, $e^{-tH}$ is a symmetric Markov semigroup, where $H=-\Delta+\frac{R}{4}$ 
in the case $\lambda_0(g_0)>0$ and $H=-\Delta+\frac{R}{4}-\frac{\min R^-}{4}$ in the general case.  
Hence the general theory and results in Chapter 2 of [D] can be applied to our situation. However, 
to obtain the precise geometric dependence of the Sobolev inequalities in 
Theorem D and Theorem $\mbox{D}^*$, one has to verify the exact geometric nature of the constants which would appear in 
the many steps of the involved (and tightly formulated) arguments in [D].   Our proofs of Theorem D and Theorem $\mbox{D}^*$  would be unclear and 
non-transparent if we go through a multitude of checking processes. Instead, we 
adapt the theory in [D] to our geometric set-up and work it out in complete, self-contained details.  
Another reason for doing so is to obtain some useful extentions of the theory 
as presented in Section 5, Appendix B and Appendix C (in particular Theorem C.5, Theorem C.6 and Theorem C.7).  On the other hand, we think that our presentation makes the theory easily accessible to the general 
audience of geometric analysis.  In particular, our presentation demonstrates in detail how the theory of the Ricci flow interacts 
with the basic theory of harmonic analysis.

Next we deduce from Theorem D a $\kappa$-noncollapsing estimate for the Ricci flow for all time which  
improves Perelman's $\kappa$-noncollapsing result [P1] for bounded time. Our estimate is independent of 
time and hence is uniform for all time. In particular, it holds both in a finite time interval and an infinite 
time interval. Moreover, our estimate provides a clear and uniform geometric dependence on the initial metric which appears to 
be optimal qualitatively. 

The $\kappa$-noncollapsing estimate below is measured relative to upper bounds of  the scalar curvature. The original 
$\kappa$-noncollapsing result of Perelman in [P1] is formulated relative to bounds for $|Rm|$. Later, 
a  $\kappa$-noncollapsing result  for bounded time measured relative to upper bounds of  the scalar curvature was 
obtained independently by Perelman (see [KL]) and 
the present author (see [Y1]).
\\

\noindent {\bf Theorem E} {\it Assume that $\lambda_0(g_0) > 0$.  Let $t \in [0, T)$. Consider the Riemannian manifold 
$(M, g)$ with $g=g(t)$. Assume $R\le \frac{1}{r^2}$ on a geodesic ball $B(x, r)$ with $r>0$. Then 
there holds 
\ba \label{noncollapse}
vol(B(x, r)) \ge \left(\frac{1}{2^{n+3}A}\right)^{\frac{n}{2}} r^n,
\ea
where $A$ is from Theorem D. In other words, the flow $g=g(t), t\in [0, T)$ is $\kappa$-noncollapsed relative to 
upper bounds of the scalar curvature on all scales. } \\

A similar  $\kappa$-noncollapsing estimate for bounded time follows from Theorem $\mbox{D}^*$, for which the 
condition $\lambda_0(g_0) \ge 0$ is not assumed.   \\

\noindent {\bf Theorem $\mbox{E}^*$} {\it Assume that $T<\infty$.  Let $L>0$ and $t \in [0, T)$. Consider the Riemannian manifold 
$(M, g)$ with $g=g(t)$. Assume $R\le \frac{1}{r^2}$ on a geodesic ball $B(x, r)$ with $0<r \le L$. Then 
there holds 
\ba \label{noncollapse*}
vol(B(x, r)) \ge \left(\frac{1}{2^{n+3}A+2BL^2}\right)^{\frac{n}{2}} r^n,
\ea
where $A$ and $B$ are from Theorem $\mbox{D}^*$. } \\

This theorem improves the previous $\kappa$-noncollapsing results for bounded time mentioned above. 
Namely it provides an explicit estimate with clear geometric dependence. Moreover, the estimate is 
uniform up to $t=0$ (under a given upper bound for $T$). (Of course, this is also the case for Theorem E.) 

We would like to point out that Theorem $\mbox{D}^*$  and Theorem $\mbox{E}^*$  lead to 
a uniform Sobolev inequality and a uniform $\kappa$-noncollapsing estimate independent of 
time for various modified Ricci flows, see the relevant results below.
 
One special consequence of Theorem E is that one can obtain smooth blow-up limits 
at $T$ both in the case $T<\infty$ and $T=\infty$, assuming that $g$ becomes singular at $T$. Previously, this was possible only at $T<\infty$
thanks to Perelman's $\kappa$-noncollapsing result. We formulate a theorem. Let $a_k$ be a sequence of positive numbers such that  $a_k \rightarrow \infty$, 
and $T_k \in (0, T)$ with $T_k \rightarrow T$.   
Consider the rescaled Ricci flows $g_k(t)=a_k g(T_k+a_k^{-1}t)$ on $M \times (-a_k T_k, 0]$. \\

\noindent {\bf Theorem F} {\it Assume that $\lambda_0(g_0) > 0$. 
Assume that $a_kT_k \rightarrow T_{\infty}<0$.  Moreover, assume that there is a sequence of points 
$x_k \in M$ with the following property. For each $L>0$ there is a positive constant 
$K$ such that $|Rm| \le K$ holds true for $g_k(t)$ on the geodesic ball of center 
$x_k$ and radius $K$, where $t \in (-a_kT_k, 0]$ is arbitrary.  Then 
a subsequence of $(M\times (-a_k T_k, 0], g_k, x_k)$ point converges smoothly to 
a pointed Ricci flow $(M_{\infty} \times (-T_{\infty}, 0], g_{\infty}, x_{\infty})$ for 
some manifold $M_{\infty}$ and $x_{\infty} \in M_{\infty}$, such that $g_{\infty}(t)$ is complete for each $t$. 
The flow $g_{\infty}$ is $\kappa$-noncollapsed relative to upper bounds of 
the scalar curvature on all scales, where  
$\kappa=2^{-\frac{n(n+3)}{2}}A^{-\frac{n}{2}}$ and $A$ is from Theorem E.
 Moreover, there holds 
for $g_{\infty}$ at all $t$
\ba
\label{flowsobolev2}
\left( \int_M |u|^{\frac{2n}{n-2}} dvol \right)^{\frac{n-2}{n}} \le 
A\int_M (|\nabla u|^2+\frac{R}{4}u^2) dvol
\ea
for all $u \in W^{1,2}(M)$.  } \\
   
By scaling invariance, the above results extend straightforwardly to  
the modified Ricci flow 
\ba \label{modified}
\frac{\partial g}{\partial t}=-2Ric+\lambda(g, t) g
\ea
with a smooth scalar function $\lambda(g, t)$ independent of $x \in M$. 
The volume-normalized 
Ricci flow 
\ba \label{volumenormalize}
\frac{\partial g}{\partial t}=-2Ric+\frac{2}{n} {\hat R} g
\ea
on a closed manifold, with $\hat R$ denoting the average scalar curvature, is  an example of the modified 
Ricci flow. The $\lambda$-normalized Ricci flow 
\ba \label{lambda}
\frac{\partial g}{\partial t}=-2Ric+\lambda g
\ea
for a constant $\lambda$ is another example. (Of course, it reduces to the Ricci flow when 
$\lambda=1$.) The normalized K\"{a}hler-Ricci flow is a special 
case of it. 

We have e.g.~the following results.  
\\

\noindent {\bf Theorem G} {\it  Theorem D and Theorem E extend to 
the modified Ricci flow.    Theorem F also extends to the case of the volume-normalized Ricci flow and 
the $\lambda$-normalized Ricci flow, when the limit flow equation is defined accordingly.} \\

Let $g=g(t)$ be a smooth solution of the modified Ricci flow 
(\ref{modified}) on $M \times [0, T)$ for some (finite or infinite) $T>0$,  with 
a given initial metric $g_0$. We set 
\ba
T^*=\int_0^T e^{-\int_0^t \lambda(g(s), s)ds} dt.
\ea

\noindent {\bf Theorem H} {\it  Assume that $T^*<\infty$.  \\
1) There are positive constants $A$ and $B$ depending only on the dimension $n$, 
a nonpositive lower bound for $R_{g_0}$, a positive lower bound for $vol_{g_0}(M)$, an upper bound for $C_S(M,g_0)$,
and an upper bound for $T^*$, such that for each $t \in [0, T)$ and 
all $u \in W^{1,2}(M)$ there holds 
\ba \label{sobolevH}
\left( \int_M |u|^{\frac{2n}{n-2}} dvol \right)^{\frac{n-2}{n}} \le 
A\int_M (|\nabla u|^2+\frac{R}{4}u^2)dvol+Be^{-\int_0^t \lambda(g(s), s)ds} \int_M u^2 dvol. 
\ea
2)  Let $L>0$ and $t \in [0, T)$. Consider the Riemannian manifold 
$(M, g)$ with $g=g(t)$. Assume $R\le \frac{1}{r^2}$ on a geodesic ball $B(x, r)$ with $0<r \le L$. Then 
there holds 
\ba \label{noncollapseH}
vol(B(x, r)) \ge \left(\frac{1}{2^{n+3}A+2Be^{-\int_0^t \lambda(g(s), s)ds}L^2}\right)^{\frac{n}{2}} r^n.
\ea
3) The conclusion of Theorem F carries over in  the case of the volume-normalized Ricci flow and 
the $\lambda$-normalized Ricci flow, when the limit flow equation is defined accordingly.
}\\

In both Theorems G and H,  the statements of Theorem F actually extend to the general case of the normalized Ricci flow under an additional assumption on the 
convergence of the rescaled $\lambda(g, t)$. We omit the statements. \\

Combining Theorems G and H with Perelman's scalar curvature estimate [ST] we obtain the 
following corollary. \\

\noindent {\bf Theorem I} {\it Let $g=g(t)$ be a smooth solution of the 
normalized K\"{a}hler-Ricci flow 
\ba
\frac{\partial g}{\partial t}=-2Ric+2\gamma g
\ea
on $M \times [0, \infty)$ with a positive first Chern class, where $\gamma$ is the positive 
constant such that the Ricci class equals $\gamma$ times the K\"{a}hler class.   (We assume that 
$M$ carries such a K\"{a}hler structure.) Then 
the Sobolev inequality (\ref{sobolevH})  holds  true with $\lambda(g(s), s)=2\lambda$. 
Moreover, there is a positive constant $L$ depending only on the initial metric $g_0=g(0)$ and the dimension $n$ such that 
the inequality (\ref{noncollapseH}) holds true for all $t \in [0, T)$ and $0<r \le L$.

If $\lambda_0(g_0)>0$, then  
the Sobolev inequality (\ref{sobolevD}) holds true for $g$.  Moreover, there is a  positive constant 
depending only on the initial metric $g_0$ and the dimension $n$ such that 
the inequality (\ref{noncollapse}) holds true for all $t \in [0, T)$ and $0<r \le L$. 
Consequently, blow-up limits of $g$ at the time infinity  
satisfy (\ref{noncollapse}) for all $r>0$ and 
the Sobolev inequality 
\ba
\label{sobolevH}
\left( \int_M |u|^{\frac{2n}{n-2}} dvol \right)^{\frac{n-2}{n}} \le 
A\int_M |\nabla u|^2 dvol
\ea
for all $u$. 
(In particular, they must be noncompact.)  } \\

Finally, we would like to mention that Theorem D and Theorem $\mbox{D}^*$ hold true for 
the Ricci flow with surgeries of Perelman [P2], with suitable modifications as stated below. 
\\

\noindent {\bf Theorem J} {\it Let $n=3$ and $g=g(t)$ be a Ricci flow with surgeries as constructed 
in [P2] on its maximal time interval $[0, T_{max})$, with suitably chosen surgery parameters.  Let $g_0=g(0)$. Let $m(t)$ denote the number of 
surgeries which are performed up to the time $t\in (0, T_{max})$. Then there holds at  each $t \in [0, T_{max})$
\ba \label{sobolevJ*}
\left( \int_M |u|^{\frac{2n}{n-2}} dvol \right)^{\frac{n-2}{n}} \le 
A(t)\int_M (|\nabla u|^2+\frac{R}{4}u^2)dvol+B(t) \int_M u^2 dvol
\ea
for all $u \in W^{1,2}(M)$, where $A(t)$ and $B(t)$ are bounded from above in terms of a nonpositive lower bound for $R_{g_0}$, a positive lower bound for $vol_{g_0}(M)$, an upper bound for $C_S(M,g_0)$, and
an upper bound for $t$.  

If $\lambda_0(g_0)>0$, then there holds at each $t \in [0, T_{max})$
\ba \label{sobolevJ}
\left( \int_M |u|^{\frac{2n}{n-2}} dvol \right)^{\frac{n-2}{n}} \le 
A(t)\int_M (|\nabla u|^2+\frac{R}{4}u^2)dvol
\ea
for all $u \in W^{1,2}(M)$, where $A(t)$ is bounded from above in terms of  a nonpositive lower bound for $R_{g_0}$, a positive lower bound for $vol_{g_0}(M)$, an upper bound for $C_S(M,g_0)$,
a positive lower bound for $\lambda_0(g_0)$, and an upper bound for $m(t)$. 

$\kappa$-noncollapsing estimates follow as before, which lead to a considerable simplification of the arguments in [P2] about preserving 
the $\kappa$-noncollapsing property after surgeries.
Similar results hold true in higher dimensions whenever similar surgeries are performed. (The constants also 
depend on the dimension $n$. 
) } \\

This result follows from Theorem D, Theorem $\mbox{D}^*$ and a general result on Sobolev inequalities 
under surgeries. The details will be presented elsewhere. In [P2], the surgery parameters are chosen 
such that several key properties of the Ricci flow are preserved after surgery. One  is the 
$\kappa$-noncollapsing property. Since the Sobolev inequalities (\ref{sobolevJ*}) and 
(\ref{sobolevJ}) are derived without using the $\kappa$-noncollapsing property, the choice of the 
surgery parameters is also simplified. The $\kappa$-noncollapsing property follows as a consequence of (\ref{sobolevJ*}) and 
(\ref{sobolevJ}).

The results in this paper (except Theorem J) extend to the dimension $n=2$. This will be presented elsewhere.
 
We would like to acknowledge that Guofang Wei first brought our attention to Zhang's paper [Z].

\sect{The Sobolev inequality}

Consider a compact  Riemannian manifold $(M, g)$ of dimension  $n\ge 3$.
Its  Poincar\'{e}-Sobolev constant (for the exponent 2) is defined to be 
\ba
C_{P,S}(M,g)=\sup\{\|u-u_M\|_{\frac{2n}{n-2}}: u \in C^1(M), \|\nabla u\|_2=1\},
\ea
where $\|u\|_p$ denotes the $L^p$ norm of $u$ with respect to $g$, i.e. 
$\|u\|_p=(\int_M |u|^pdvol)^{1/p}$ ($dvol=dvol_g$).  
In other words, $C_{P,S}(M,g)$ is the smallest number  such that the 
Poincare-Sobolev inequality 
\ba 
\|u-u_M\|_{\frac{2n}{n-2}} \le C_{P,S}(M,g) \|\nabla u\|_2
\ea
holds true for all $u \in C^1(M)$ (or all $u \in W^{1,2}(M)$).
The Sobolev constant of $(M,g)$ (for the exponent 2) is defined to be 
\ba
C_S(M,g)=\sup\{\|u\|_{\frac{2n}{n-2}}-\frac{1}{vol(M)^{\frac{1}{n}}}\|u\|_2:
u \in C^1(M), \|\nabla u\|_2=1.
\ea
In other words, $C_S(M,g)$ is the smallest number such that the  
inequality 
\ba
\|u\|_{\frac{2n}{n-2}} \le C_{S}(M,g)\|\nabla u\|_2+\frac{1}{vol(M)^{\frac{1}{n}}}\|u\|_2
\ea
holds true for all $u\in W^{1,2}(M)$. \\

\noindent {\bf Definition} We define the {\it modified Sobolev constant} 
$\tilde C_S(M,g)$ to be $\max\{C_S(M,g), 1\}$. \\

The H\"{o}lder inequality leads to the following basic fact.

\begin{lem} There holds for all $u\in W^{1,2}(M)$
\ba
\|u\|_{\frac{2n}{n-2}} \le C_{P,S}(M,g)\|\nabla u\|_2+\frac{1}{vol(M)^{\frac{1}{n}}}\|u\|_2.
\ea
In other words, there holds $C_S(M,g) \le C_{P,S}(M,g)$.
\end{lem}

Another basic constant, the Neumann isoperimetric constant of 
$(M, g)$,  is defined to 
be 
\ba
C_{N,  I}(M,g)=\sup\{\frac{vol(\Omega)^{\frac{n-1}{n}}}{A(\partial \Omega)}: 
\Omega  \subset M \mbox{ is a } 
C^1 \mbox{ domain }, vol(\Omega) \le \frac{1}{2}  vol(M)\},
\ea
where $A(\partial \Omega)$ denotes the $n-1$-dimensional volume of $\partial \Omega$.

\begin{lem} \label{poincare} There holds for all $u \in  W^{1,2}(M)$
\ba \label{poincare2}
\|u-u_M\|_{\frac{2n}{n-2}} \le  2(1+\sqrt{2}) \frac{n-1}{n-2} C_{N,I}(M,  g)
\|\nabla u\|_2.
\ea
In other words, there holds $C_{P,S}(M,g) \le 2(1+\sqrt{2})\frac{n-1}{n-2} C_{N,I}(M,  g)$.
\end{lem} 

For the proof see [Y3]. The following estimate of the Neumann isoperimetric constant follows from   
S.~Gallot's estimate in [Ga2]. We define the diamater rescaled Ricci curvature 
$\hat Ric(v,v)$ of a unit tangent vector $v$ to be $diam(M)^2 Ric(v,v)$, and set 
$\kappa_{\hat Ric}=\min_v\{\hat Ric(v,v)\}$. Then we set
$\hat \kappa_{\hat Ric}=|\min\{\kappa_{\hat Ric}, -1\}|$. 
We also define the diameter rescaled volume $\hat vol(M)$ to be $vol(M)diam(M)^{-n}.$
 
\begin{theo}  \label{gallot} There holds
\ba  \label{gallot}
C_{N,I}(g, M) \le C(n, \hat \kappa_{\hat  Ric})\hat vol(M)^{-\frac{1}{n}},
\ea
where $C(n, \hat \kappa_{\hat Ric})$ is a positive constant  depending only 
on $n$ and
$\hat \kappa_{\hat Ric}$.  
\end{theo} 

Note that $\hat \kappa_{\hat Ric}$ can be replaced by a certain integral 
lower bound of the Ricci curvature, see [Ga1].

\sect{The Logarithmic Sobolev inequalities on a Riemannian Manifold} 

The various versions of the logarithmic Sobolev inequality on 
the Euclidean space as presented in Appendix A allow suitable extentions to Riemannian manifolds. 
We formulate a log gradient version and a straight version, cf. Appendix A. As in the last section, let $(M, g)$ be 
a compact Riemannian manifold of dimension $n$. 

\begin{theo} \label{logtheorem1} There holds 
\ba \label{RLS1}
\int_M u^2 \ln u^2 dvol \le n \ln \left(C_S(M,g) \|\nabla u\|_2+ \frac{1}{vol_g(M)^{\frac{1}{n}}}\right),
\ea
provided that $u \in W^{1,2}(M)$ and $\|u\|_2=1$. 
\end{theo}

\Pf  Set $q=\frac{2n}{n-2}$.  Since $\ln $ is concave and 
$\int_M u^2 dvol =1$, we have by Jensen's inequality
\ba
\ln \int_M u^q dvol = \ln \int_M u^2 \cdot u^{q-2} dvol \ge \int_M u^2 \ln u^{q-2}.
\ea
It follows that 
\ba
\int_M u^2 \ln u &\le& \frac{1}{q-2} \ln \int_M u^q  dvol \nonumber \\
&=& \frac{q}{q-2} \ln \|u\|_q \nonumber \\
&\le& \frac{n}{2} \ln  \left(C_S(M,g) \|\nabla u\|_2+ \frac{1}{vol_g(M)^{\frac{1}{n}}}\|u\|_2\right).
\ea 
\qed

\begin{lem} \label{loglemma1}  There holds 
\ba \label{log1}
\ln (x+B) \le \alpha x+\alpha B-1-\ln \alpha
\ea
for all $B \ge 0, \alpha>0$ and $x >-B$.
\end{lem}
\Pf Consider the function $y=\ln (x+B)-\alpha x$ for $x>-B$. Since $y \rightarrow -\infty$ as 
$x \rightarrow -B$ or $x \rightarrow \infty$, it achieves its maximum somewhere. We have 
\ba
y'=\frac{1}{x+B}-\alpha.
\ea
Hence the maximum point is $x_0=\frac{1}{\alpha}-B$. It follows that the maximum of $y$ is 
$y(x_0)=\alpha B-1-\ln \alpha$. 
\qed\\

\begin{theo}  \label{logtheorem2} \label{RLS2-th} For each $\alpha>0$ and all $u \in W^{1,2}(M)$ with $\|u\|_2=1$ there holds
\ba \label{RLS2}
\int_M u^2 \ln u^2  \le \frac{n\alpha C_S(M,g)^2}{2} \int_M |\nabla u|^2 -\frac{n}{2}\ln \alpha+ \frac{n}{2}(\ln 2+\alpha vol_g(M)^{-\frac{2}{n}}-1)
\ea
and
\ba \label{RLS3}
\int_M u^2 \ln u^2  &\le& \frac{n\alpha C_S(M,g)^2}{2} \int_M (|\nabla u|^2+\frac{R}{4}u^2) -\frac{n}{2}\ln \alpha
\nonumber \\
&& +\frac{n\alpha}{2}(vol_g(M)^{-\frac{2}{n}}-\frac{\min R^-}{4}C_S(M,g)^2)+\frac{n}{2}(\ln 2-1).
\ea
(The notation of the volume is omitted.) 
\end{theo}
\Pf By (\ref{RLS1}) we have for $u \in W^{1,2}(M)$ with $\|u\|_2=1$
\ba \label{middle1}
\int_M u^2 \ln u^2 &\le& \frac{n}{2} \ln \left(C_S(M,g) \|\nabla u\|_2+ \frac{1}{vol_g(M)^{\frac{1}{n}}}\right)^2 
\nonumber \\ &\le&  \frac{n}{2} \ln 2 + \frac{n}{2} \ln \left(C_S(M,g)^2 \int_M |\nabla u|^2+ \frac{1}{vol_g(M)^{\frac{2}{n}}}\right).
\ea
Applying  Lemma  \ref{loglemma1} with $x=C_S(M,g)^2 \int_M |\nabla u|^2$ and $B=1$ 
we then arrive at 
(\ref{RLS2}). The inequality (\ref{RLS3}) follows from (\ref{RLS2}). \qed \\

\begin{lem} Let $A>0, B>0$ and $\gamma>0$ such that $A\ge \frac{1}{\gamma+B}$. Then we have 
\ba \label{log2}
\ln(x+B)\le Ax-\ln A+\ln (\gamma+B)-\ln \gamma-1
\ea
for all $x \ge \gamma$.
\end{lem}
\Pf First consider the function $y=\ln t-\gamma t$ for $t>0$. Since $y\rightarrow  -\infty$ as 
$t \rightarrow 0$ or $t \rightarrow \infty$, $y$ achieves its maxmum somewhere. We have 
$y'=\frac{1}{t}-\gamma$. Hence the maximum is achieved at $\frac{1}{\gamma}$. It follows that 
the maximum is $y(\frac{1}{\gamma})=-\ln \gamma -1$. We infer 
\ba \label{log-A}
\ln A-\gamma A \le -\ln \gamma -1.
\ea

Next we consider the function $y=\ln (x+B)-Ax+\ln A$ for $x \ge \gamma$. By (\ref{log-A}) we have 
$y(\gamma)=\ln (\gamma+B)-A\gamma+\ln A \le \ln (\gamma+B)-\ln \gamma -1$. On the other hand, 
we have $y'=\frac{1}{x+B}-A \le \frac{1}{\gamma+B}-A \le 0$. We arrive at (\ref{log2}). 
\qed\\

\begin{theo} \label{logtheorem3} Assume that the first eigenvalue $\lambda_0=\lambda_0(g)$ of the operator 
$-\Delta+\frac{R}{4}$ is positive.  For each $A \ge \delta_0$ and all 
$u \in W^{1,2}(M)$ with $\|u\|_2=1$ there holds 
\ba \label{RLS4}
\int_M u^2 \ln u^2 \le \frac{nAC_S^2}{2}\int_M(|\nabla u|^2+\frac{R}{4}u^2)
-\frac{n}{2}\ln A+\frac{n}{2}\ln 2+\sigma_0,
\ea
where
\ba \label{delta-0}
\delta_0=\delta_0(g)=
(\lambda_0 C_S^2+\frac{1}{vol_g(M)^{\frac{2}{n}}}
-C_S^2 \frac{\min R^-}{4})^{-1},
\ea 
\ba \label{sigma-0}
\sigma_0=\sigma_0(g)==\frac{n}{2}\left[\ln (\lambda_0 C_S^2+\frac{1}{vol_g(M)^{\frac{2}{n}}}
-C_S^2 \frac{\min R^-}{4})-\ln (\lambda_0 C_S^2)-1\right], 
\ea
and $C_S=C_S(M,g)$.
\end{theo}
\Pf Arguing as in the proof of Theorem \ref{RLS2-th} we deduce for $u \in W^{1,2}(M)$ with $\|u\|_2=1$
\ba
&&\int_M u^2 \ln u^2 \le \frac{n}{2}\ln 2+ \frac{n}{2} \ln(C_S^2\int_M |\nabla u|^2 +\frac{1}{vol_g(M)^{\frac{2}{n}}})
\nonumber \\
&\le& \frac{n}{2}\ln 2+ \frac{n}{2} \ln \left[C_S^2\int_M (|\nabla u|^2+\frac{R}{4}u^2) +\frac{1}{vol_g(M)^{\frac{2}{n}}}
-C_S^2 \frac{\min R^-}{4}\right].
\ea
Applying (\ref{log2}) with $\gamma=\lambda_0 C_S^2$, $B=\frac{1}{vol_g(M)^{\frac{2}{n}}}
-C_S^2 \frac{\min R^-}{4}$ and $x=C_S^2\int_M (|\nabla u|^2+\frac{R}{4}u^2)$ we then 
arrive at (\ref{RLS4}) for each $A \ge (\gamma+B)^{-1}$. \qed

\sect{The logarithmic Sobolev inequality along the Ricci flow }

Let $M$ be a compact manifold of dimension $n$. Consider Perelman's 
entropy functional 
\ba \label{entropy}
{\mathcal W}(g, f, \tau)=\int_M \left[ \tau(R+|\nabla f|^2)+f-n\right] \frac{e^{-f}}{(4\pi\tau)^{\frac{n}{2}}} dvol,
\ea
where $\tau$ is a positive number, $g$ is a Riemannian metric on $M$,  and $f\in C^{\infty}(M)$ satisfies 
\ba \label{vol-1}
\int_M \frac{e^{-f}}{(4\pi\tau)^{\frac{n}{2}}} dvol=1.
\ea
All goemetric quantities in (\ref{entropy}) and (\ref{vol-1}) are associated with $g$. 
To relate to the idea of logarithmic Sobolev inequalities we make a change of 
variable 
\ba
u=\frac{e^{-\frac{f}{2}}}{(4\pi\tau)^{\frac{n}{4}}}.
\ea
Then (\ref{vol-1}) leads to 
\ba \label{vol-2}
\int_M u^2 dvol=1
\ea
 and we have
\ba
{\mathcal W}(g, f, \tau)={\mathcal W}^*(g, u, \tau) -\frac{n}{2}\ln \tau-\frac{n}{2}\ln(4\pi)-n
\ea
where 
\ba
{\mathcal W}^*(g, u,\tau)=\int_M \left[\tau(4|\nabla u|^2+R u^2)-u^2 \ln u^2 \right] dvol.
\ea
We define $\mu^*(g, \tau)$ to be the infimum of ${\mathcal W}^*(g, u, \tau)$ over all 
$u$ satisfying (\ref{vol-2}).
 
Next let 
$g=g(t)$ be a smooth solution of the Ricci flow 
\ba
\frac{\partial g}{\partial t}=-2Ric
\ea
on $M\times [0, T)$ for some (finite or infinite) $T>0$.
Let $0<t^*<T$ and $\sigma>0$.  We set $T^*=t^*+\sigma$ and 
$\tau=\tau(t)=T^*-t$ for $0\le t \le t^*$. 
Consider a solution $f=f(t)$ of the equation
\ba \label{nonlinearconjugate}
\frac{\partial f}{\partial t}=-\Delta f+|\nabla f|^2-R+\frac{n}{2\tau}
\ea
on $[0, t^*]$ with a given terminal value at $t=t^*$ (i.e. $\tau=\sigma$)
satisfying (\ref{vol-1}) with $g=g(t^*)$. Then (\ref{vol-1}) holds true for 
$f=f(t), g=g(t)$ and all $t \in [0, t^*]$.
Perelman's 
monotonicity formula says 
\ba
\frac{d {\mathcal W}}{d t}=2\tau \int_M |Ric+\nabla^2 f-\frac{1}{2\tau} g|^2 \frac{e^{-f}}{(4\pi\tau)^{\frac{n}{2}}} dvol
\ge 0,
\ea
where ${\mathcal W}={\mathcal W}(g(t), f(t), \tau(t))$. Consequently, 
\ba
\frac{d}{dt} {\mathcal W}^*(g,u,\tau) \ge \frac{n}{2}\frac{d}{dt} \ln \tau,
\ea
where $g=g(t), \tau=\tau(t)$ and 
\ba
u=u(t)=
\frac{e^{-f(t)/2}}{(4\pi\tau(t))^{\frac{n}{4}}},
\ea
which satisfies the equation
\ba
\frac{\partial u}{\partial t}=-\Delta u+\frac{|\nabla u|^2}{u}+\frac{R}{2}u.
\ea

It follows that 
\ba
\mu^*(g(t_1), \tau(t_1)) \le \mu^*(g(t_2, \tau(t_2))+\frac{n}{2} \ln \frac{\tau_1}{\tau_2},
\ea
for $t_1<t_2$, where $\tau_1=\tau(t_1)$ and $\tau_2=\tau(t_2)$. Choosing $t_1=0$ and 
$t_2=t^*$ we then arrive at 
\ba \label{mu-1}
\mu^*(g(0), t^*+\sigma) \le \mu^*(g(t^*), \sigma)+\frac{n}{2} \ln \frac{t^*+\sigma}{\sigma}.
\ea
Since $0<t^*<T$ is arbitrary, we can rewrite (\ref{mu-1}) as follows
\ba \label{mu-monotone}
\mu^*(g(t), \sigma) \ge \mu^*(g(0), t+\sigma)+ \frac{n}{2}\ln \frac{\sigma}{t+\sigma}
\ea
for all $ t \in [0, T)$ and $\sigma>0$ (the case $t=0$ is trivial). 
\\

We'll also need the following elementary lemma. 

\begin{lem} \label{stronglemma} Let $a>0$ and $b$ be constants. Then the minimum of the 
function $y=a \sigma-\frac{n}{2} \ln \sigma +b$ for $\sigma>0$ is 
$\frac{n}{2} \ln(\alpha a)$, where 
\ba \label{strong}
\alpha=\frac{2e}{n} e^{\frac{2b}{n}}.
\ea
\end{lem}  
\Pf Since $y \rightarrow \infty$ as $t\rightarrow 0$ or $t \rightarrow \infty$, it 
achieves its minimum somewhere. We have $y'=a-\frac{n}{2\sigma}$, whence the 
minimum is achieved at $\sigma=\frac{n}{2a}$. Then the minimum equals 
$y(\frac{n}{2a})$, which leads to  the desired conclusion.  \qed \\

\noindent {\bf Proof of Theorem A} \\

We apply  Theorem \ref{logtheorem2} with $g=g_0$ to estimate $\mu^*(g_0, t+\sigma)$. Consider  $u \in W^{1,2}(M)$ with 
$\|u\|_2=1$.  We choose 
\ba
\alpha=\frac{8(t+\sigma)}{n\tilde C_S(M,g_0)^2}
\ea
in (\ref{RLS2}) and deduce
\ba
\int_M u^2 \ln u^2 &\le& 4(t+\sigma)\int_M |\nabla u|^2 -\frac{n}{2}\ln \frac{8(t+\sigma)}{n\tilde C_S^2}
+\frac{n}{2} \cdot \frac{8(t+\sigma)}{n\tilde C_S^2 vol_{g_0}(M)^{\frac{2}{n}}}+\frac{n}{2}(\ln 2-1)
\nonumber \\ &\le& (t+\sigma) \int_M (4|\nabla u|^2+Ru^2) 
+(t+\sigma)(\frac{4}{n\tilde C_S^2 vol_{g_0}(M)^{\frac{2}{n}}}-\min_{t=0} R) 
\nonumber \\ &&-\frac{n}{2} \ln (t+\sigma)+\frac{n}{2}(2\ln \tilde C_S+\ln n -2\ln2 -1),
\ea
where $\tilde C_S=\tilde C_S(M,g_0)$. 
It follows that 
\ba
\mu^*(g(0), t+\sigma) &\ge& \frac{n}{2} \ln (t+\sigma) -(t+\sigma)(\frac{4}{n\tilde C_S^2 vol_{g_0}(M)^{\frac{2}{n}}}-\min_{t=0} R) 
\nonumber \\ &&-\frac{n}{2}(2\ln \tilde C_S+\ln n -2\ln2 -1).
\ea
Combining this with (\ref{mu-monotone}) leads to 
\ba 
\mu^*(g(t), \sigma) &\ge& \frac{n}{2} \ln \sigma -(t+\sigma)(\frac{4}{n\tilde C_S^2 vol_g(M)^{\frac{2}{n}}}-\min_{t=0} R) 
\nonumber \\ &&-\frac{n}{2}(2\ln \tilde C_S+\ln n -2\ln2 -1),
\ea
or 
\ba 
\mu^*(g(t), \frac{\sigma}{4}) &\ge& \frac{n}{2} \ln \sigma -(t+\frac{\sigma}{4})(\frac{4}{n\tilde C_S^2 vol_g(M)^{\frac{2}{n}}}-\min_{t=0} R) 
\nonumber \\ &&-\frac{n}{2}(2\ln \tilde C_S+\ln n-1),
\ea
which is equivalent to (\ref{sobolevA}). 

To see (\ref{strongA}) we apply Lemma \ref{stronglemma} to 
(\ref{sobolevA}) with 
$a=\int_M(|\nabla u|^2+\frac{R}{4}u^2)dvol+\frac{A_1}{4}$ and 
$b=A_1t+A_2$. \qed \\

\noindent {\bf Proof of Theorem B} \\

This is similar to the proof of Theorem A. We apply Theorem \ref{logtheorem3} with $g=g_0$  
to estimate $\mu^*(g_0, t+\sigma)$. Assume $t+\sigma \ge \frac{n}{8}C_S(M,g_0)^2\delta_0(g_0)$.
We set
\ba
A=\frac{8(t+\sigma)}{nC_S(M,g_0)^2}.
\ea
Then there holds $A \ge \delta_0(g_0)$. Using this $A$ in (\ref{RLS4}) we deduce for   $u \in W^{1,2}(M)$ with 
$\|u\|_2=1$
\ba
\int_M u^2 \ln u^2 &\le& 4(t+\sigma)\int_M (|\nabla u|^2+\frac{R}{4}u^2)
-\frac{n}{2}\ln (t+\sigma)\nonumber \\
&&+\frac{n}{2}(2\ln C_S(M,g_0)+\ln n-2\ln 2)+\sigma_0(g_0).
\ea
It follows that 
\ba
\mu^*(g_0, t+\sigma) \ge \frac{n}{2}\ln (t+\sigma)-\frac{n}{2}(2\ln C_S(M,g_0)+\ln n-2\ln 2)-\sigma_0(g_0).
\ea
Combining this with (\ref{mu-monotone}) yields 
\ba
\mu^*(g(t), \sigma)\ge \frac{n}{2} \ln \sigma -\frac{n}{2}(2\ln C_S(M,g_0)+\ln n-2\ln 2)-\sigma_0(g_0).
\ea
Replacing $\sigma$ by $\frac{\sigma}{4}$ we then arrive at (\ref{sobolevB}).

To see (\ref{strongB}), we apply Lemma \ref{stronglemma} to (\ref{sobolevB}) with 
$a=\int_M (|\nabla u|^2+\frac{R}{4}u^2)$ and $b=\frac{n}{2}\ln n+n \ln C_S(M, g_0)+\sigma_0(g_0)$. 
Note that by the maximum principle and the evolution equation of the scalar curvature associated with the 
Ricci flow, $\min R$ is nondecreasing, which implies that $a>0$.  \qed \\

Note that the proofs of Theorem A and Theorem B lead to the following general result.
Indeed, Theorem A and Theorem B follows from it.  

\begin{theo} \label{riccilog} Let $g=g(t)$ be a smooth solution of the Ricci flow on 
$M \times [0, T)$ for some (finite or infinite) $T>0$. Let $h(\sigma)$ be 
a scalar function for $\sigma>0$.  Assume that 
the initial metric $g_0=g(0)$ satisfies the logarithmic Sobolev inequality 
\ba
\int_M u^2 \ln u^2 dvol \le \sigma \int_M (|\nabla u|^2 + \frac{R}{4}u^2)dvol
+h(\sigma)
\ea
for each $\sigma>0$ and all $u\in W^{1,2}(M)$ with $\int_M u^2 dvol=1$. Then 
there holds at each $t \in [0, T)$
\ba
\int_M u^2 \ln u^2 dvol \le \sigma \int_M (|\nabla u|^2 + \frac{R}{4}u^2)dvol
+h(4(t+\sigma))
\ea
for each $\sigma>0$ and all $u\in W^{1,2}(M)$ with $\int_M u^2 dvol=1$.
\end{theo}

\noindent {\bf Proof of Theorem C} \\

We first consider the case $\lambda_0(g_0)>0$.  Let $t\in [0, T)$ and $\sigma>0$. If $\sigma < \frac{n}{8}C_S(M, g_0)^2 \delta_0(g_0)$, we apply Theorem A. 
Otherwise, we apply Theorem B. Then we arrive at (\ref{sobolevC}).  To see (\ref{strongC}), we note that 
by [P1] the eigenvalue $\lambda_0(g(t))$ is nondecreasing. Hence 
$\lambda_0(g(t))>0$ for all $t$, which implies that $\int_M (|\nabla u|^2+\frac{R}{4}u^2)>0$ for all $t$.  Hence
we can apply Lemma \ref{stronglemma} to (\ref{sobolevC}) with $a=\int_M (|\nabla u|^2+\frac{R}{4}u^2)$ and $b=C$
to arrive at the desired inequality. 

Next we consider the case $\lambda_0(g_0)=0$. Consider $t_0=\min\{\frac{T}{2}, 1\}$. There are two cases 
to consider. The first case is
$\lambda_0(g(t_0)>0$.  In this case we first apply Theorem A to obtain (\ref{sobolevC}) 
for $0 \le t \le t_0$.  Then we repeat the arguments in the proof of Theorem B, using the 
derived logarithmic Sobolev inequality at $t_0$ instead of Theorem \ref{?}, with $t_0$ as the initial time.  
Then we arrive at (\ref{sobolevC}) for $t\ge t_0$.  The inequality (\ref{strongC}) in the same way as 
before.
The second case is $\lambda_0(g(t_0))=0$. In this case, $g=g(t)$ is a steady gradient Ricci soliton 
on $[0, t_0]$. Hence it is also a steady gradient Ricci soliton on $[0, T)$.  It follows that 
the logarithmic Sobolev inequalities at $t=0$ provided by Theorem  and Theorem remain unchanged 
for $t>0$.  \qed \\

\noindent {\bf Proof of Corollary to Theorem C} \\

Choosing $u=vol_{g(t)}(M)^{-\frac{1}{2}}$ in (\ref{sobolevC}) we infer 
\ba
\ln \frac{1}{vol_{g(t)}(M)}\le \frac{\sigma}{4}\hat R(t)-\frac{n}{2}\ln \sigma +C.
\ea
If $\hat R(t) \le 0$ we choose $\sigma=1$ to arrive at (\ref{volume3}). If $\hat R(t)>0$, we choose 
$\sigma=\hat R(t)^{-1}$ to arrive at (\ref{volume4}). \qed

\sect{The Sobolev inequality along the Ricci flow}

We first present a general result which converts a logarithmic 
Sobolev inequality to a Sobolev inequality. It follows straightforwardly from more general results in [D].
Consider a compact Riemannian 
manifold $(M,g)$ of dimension $n \ge 1$. Let $\Psi \in L^{\infty}(M)$, which we call a potential 
function.
We set $H=-\Delta+\Psi$. Its associated quadratic form is 
\ba
Q(u)=\int_M (|\nabla u|^2 +\Psi u^2) dvol,
\ea
where $u \in W^{1,2}(M)$.  We also use $Q$ to denote the corresponding bilinear form, i.e. 
\ba
Q(u, v)=\int_M (\nabla u \cdot \nabla v +\Psi u v)dvol.
\ea
 Consider the operator  $e^{-tH}$ associated with  
$H$.  It is characterized by the property that 
for $u_0 \in L^2(M)$, $u=e^{-tH}u_0$ satisfies the heat equation
\ba
\frac{\partial u}{\partial t}=-Hu
\ea
for $t>0$ and the initial condition $u(0)=u_0$.   We have the spectral formula 
\ba \label{spectral}
e^{-tH}u=\sum e^{-\lambda_i t} \phi_i <u, \phi_i>_2,
\ea
for $u \in L^2(M)$, where $\{\phi_i\}$ is a complete 
set of $L^2$-orthonormal eigenfunctions of $H$ and $\lambda_1\le \lambda_2 \le \cdot \cdot \cdot$ are 
the corresponding eigenvalues.  Since $\lambda_i \rightarrow \infty$, $e^{-tH}: L^2(M) 
\rightarrow L^2(M)$ is a bounded operator.  On the other hand, there holds 
\ba
e^{-tH}u=\int_M K(\cdot, y, t)u dvol_y,
\ea
where $K(x,y,t)$ denotes the heat kernel of $H$. 

\begin{lem} \label{extensionlemma}  The extension of $e^{-tH}$ for $t>0$ to $L^1(M)$ by the spectral formula 
(\ref{spectral}) defines a bounded linear operator $e^{-tH}: 
L^1(M) \rightarrow W^{2, p}(M)$ for each $ 0<p<\infty$. 
\end{lem}
\Pf By elliptic regularity we have $\phi_i \in W^{2, p}(M)$ for each $i$ and 
$0<p<\infty$.  The elliptic $W^{2,p}$ estimates and Sobolev embedding lead to 
$\|\phi_i\|_{2, p} \le c_p(|\lambda_i|+1)^{m_n} $ for some $c_p>0$ independent of $i$ and 
a natural number $m_n$ depending only on $n$.  The Sobolev embedding then implies 
$\|\phi_i\|_{\infty} \le c(\lambda_i|+1)^{m_n} $ for some $c>0$ independent of $i$. 
Now we have for $u \in L^1(M)$ 
\ba
\sum_{i\ge 1} e^{-\lambda_i t} |<u, \phi_i>_2 | \cdot \|\phi_i\|_{2, p} \le 
\left(\sum_{i \ge 1} e^{-\lambda_i t} \|\phi_i\|_{\infty} \|\phi_i\|_{2, p}\right) \|u\|_1.
\ea
By the above estimates, the last series converges.  The desired conclusion follows. \qed \\

\begin{lem} \label{contractionlemma} Assume $\Psi \ge 0$. Then $e^{-tH}$ for $t>0$ is a contraction on $L^p(M)$ for each 
$1\le p \le \infty$, i.e. 
\ba
\|e^{-tH}u\|_p \le \|u\|_p
\ea
for all $u \in L^p(M)$.  It is also a contraction on $W^{1,2}(M)$ with respect to the norm 
$Q(u)^{\frac{1}{2}}$ (if $Q>0$, i.e. $\lambda_1>0$) or 
the norm $(Q(u)+\int_M u^2 dvol)^{\frac{1}{2}}$ (if $\lambda_1=0$).  Moreover, it is positivity preserving, i.e. 
$e^{-tH}u \ge 0$ if $u \ge 0$ and $u \in L^2(M)$.  
\end{lem}
\Pf The maximum principle implies that $e^{-tH}$ is a contraction on $L^{\infty}(M)$ for 
$t>0$. For $t>0$ and $u \in L^1(M)$ we set $\phi=sgn(e^{-tH}u)$, i.e. $\phi=1$ where $e^{-tH}u \ge 0$ and $\phi=-1$ where
$e^{-tH}u<0$. There holds 
\ba
\|e^{-tH}u\|_1 &=&
 \int_M \phi e^{-tH}u  dvol = \int_M u e^{-tH} \phi dvol 
 \nonumber \\  &\le& \|e^{-tH}\phi\|_{\infty} \|u\|_1
\le \|\phi\|_{\infty} \|u\|_1=\|u\|_1. 
\ea
Hence $e^{-tH}$ is a contraction on $L^1(M)$.  By the Riesz-Thorin interpolation theorem (see Appendix C),  
$e^{-tH}$ is a contraction on $L^p(M)$ for each $1< p < \infty$. 

The contraction property of $e^{-tH}$ on 
$W^{1,2}(M)$ follows from the spectral formula (\ref{spectral}) because $\lambda_1 \ge 0$. 
(The contraction property of $e^{-tH}$ on $L^2(M)$ also follows from  
(\ref{spectral}).)  Finally, the positivity preserving property of 
$e^{-tH}$ is a consequence of the maximum principle. \qed \\

\begin{theo} \label{D1} Let $0<\sigma^*\le \infty$. Assume that for each $0<\sigma<\sigma^*$ the logarithmic Sobolev 
inequality 
\ba \label{Dlog1}
\int_M u^2 \ln u^2 dvol \le \sigma Q(u)+ \beta(\sigma)
\ea
holds true for all $u \in W^{1,2}(M)$ with $\|u\|_2=1$, where 
$\beta$ is a non-increasing continuous function.  Assume that 
\ba
\tau(t)=\frac{1}{2t}\int^t_0 \beta(\sigma)d\sigma
\ea
is finite for all $0<t < \sigma^*$. Then there holds 
\ba \label{heat1}
\|e^{-tH}u\|_{\infty} \le e^{\tau(t)-\frac{3t}{4}\inf \Psi^-} \|u\|_2
\ea
for each $0<t< \frac{1}{4}\sigma^*$ and all $u \in L^2(M)$. 
There also holds 
\ba \label{heat2}
\|e^{-tH}u\|_{\infty} \le e^{2\tau(\frac{t}{2})-\frac{3t}{4} \inf \Psi^-} \|u\|_1
\ea
for each $0<t< \frac{1}{4}\sigma^*$ and all $u \in L^1(M)$. 
 \end{theo}

The proof of this theorem is presented in Appendix B. Note that (\ref{heat2}) is equivalent to an upper bound for the heat kernel.   The  nonincreasing condition on $\beta$ 
can easily be removed (the function $\tau(t)$ needs to be slightly modified).

\begin{theo} \label{D2} 1) Assume $\Psi \ge 0$.  Let $\mu>2$ and $c>0$. Assume that the inequality 
\ba
\|e^{-tH}u\|_{\infty} \le c t^{-\frac{\mu}{4}} \|u\|_2
\ea 
holds true for each $t>0$ and all $u\in L^2(M)$. Then   
the Sobolev inequality 
\ba \label{Dsobolev2}
\|u\|^2_{\frac{2\mu}{\mu-2}} \le C(\mu, c) Q(u)
\ea
holds true for all $u \in W^{1,2}(M)$, where the positive constant 
$C(\mu, c)$ can be bounded from above in terms of 
upper bounds for $c, \mu$ and $\frac{1}{\mu-2}$. \\
2) Let $\mu>2$ and $c>0$.  Assume that the inequality 
\ba 
\|e^{-tH}u\|_{\infty} \le c_1 t^{-\frac{\mu}{4}} \|u\|_2
\ea 
holds true for each $0<t<1$ and all $u\in L^2(M)$. Then   
the Sobolev inequality 
\ba \label{Dsobolev1}
\|u\|^2_{\frac{2\mu}{\mu-2}} \le C(\mu, c) (Q(u)+(1-\inf \Psi^-)\|u\|_2^2)
\ea
holds true for all $u \in W^{1,2}(M)$, where $C(\mu, c)$ has the same property as 
the above $C(\mu, c)$.
\end{theo} 

The proof of this theorem is presented in Appendix C.
Combining Theorem \ref{D1} and \ref{D2} we arrive at the following result. 

\begin{theo} \label{D3}
Let $0<\sigma^*<\infty$. Assume that for each $0<\sigma< \sigma^*$ the logarithmic Sobolev 
inequality 
\ba \label{Dlog33}
\int_M u^2 \ln u^2 dvol \le \sigma Q(u)-\frac{\mu}{2} \ln \sigma+C
\ea
holds true for all $u \in W^{1,2}(M)$ with $\|u\|_2=1$, where 
$\mu$ and $c$ are constants such that $\mu>2$.    Then we have the Sobolev inequality
\ba \label{Dsobolev3}
\|u\|^2_{\frac{2\mu}{\mu-2}} \le  \left(\frac{\sigma^*}{4}\right)^{1-\frac{n}{\mu}} C(\bar C, \mu)
\left(Q(u)+\frac{4-\sigma^* \min \Psi^-}{\sigma^*} \|u\|_2^2\right)
\ea
for all $u \in W^{1,2}(M)$, where $C(\bar C, \mu)$ is from Theorem \ref{D2} and 
$\bar C$ is defined in (\ref{barC}) below. 
\end{theo}
\Pf For $\lambda>0$ we consider the metric $\bar g=\lambda^{-2}g$ and the 
potential function $\bar \Psi=\lambda^2 \Psi$. Let $\bar H=-\Delta_{\bar g}+
\bar \Psi$ and $\bar Q$ the associated quadratic form. It follows from 
(\ref{Dlog33}) that
\ba \label{Dlog34}
\int_M u^2 \ln u^2 dvol_{\bar g} 
\le \sigma \bar Q(u)-\frac{\mu}{2} \ln \sigma +(n-\mu)\ln \lambda +C
\ea 
for $0<\sigma < \lambda^{-2} \sigma^*$ and $u \in W^{1,2}(M)$ with 
$\|u\|_2=1$.  Choosing $\lambda=\frac{1}{2}\sqrt{\sigma^*}$ we obtain 
\ba \label{Dlog35}
\int_M u^2 \ln u^2 dvol_{\bar g} 
\le \sigma \bar Q(u)-\frac{\mu}{2} \ln \sigma +\frac{n-\mu}{2}(\ln \sigma^*
-2\ln 2) +C
\ea 
for each $0<\sigma <4$. By Theorem \ref{D1} we have for each $0<t<1$
and $u \in L^2(M)$  
\ba
\|e^{-tH}u\|_{\infty} \le \bar C
t^{-\frac{\mu}{4}} \|u\|_{2, \bar g},
\ea
where 
\ba \label{barC}
\bar C=2^{\frac{\mu-n}{2}}(\sigma^*)^{\frac{n-\mu}{4}}e^{\frac{\mu}{4}-\frac{3\sigma^*}{16}\min \Psi^-+\frac{1}{2}C}. 
\ea
Applying  Theorem \ref{D2} and converting back to $g$ we then arrive at (\ref{Dsobolev3}). \qed \\

\noindent {\bf Proof of Theorem D} \\

Applying Theorem C and Theorem \ref{D3} with $\Psi=\frac{R}{4}$,$\mu=n$ and 
$\sigma^*=4$ we deduce 
\ba
\|u\|^2_{\frac{2n}{n-2}} \le c \left(\int_M (|\nabla u|^2+\frac{R}{4}u^2)dvol
+(1-\frac{\min_{t} R^-}{4})\int_M u^2 dvol \right),
\ea
where $c=c(C, -\min_t R^-)$. By the maximum principle, we have 
$\min_t R^- \ge \min_{t=0}  R^-$. Hence we arrive at 
\ba
\|u\|^2_{\frac{2n}{n-2}} \le c \left(\int_M (|\nabla u|^2+\frac{R}{4}u^2)dvol
+(1-\frac{\min_{0} R^-}{4})\int_M u^2 dvol \right)
\ea
with $c=c(C, -\min_0 R^-)$.
 Since $\lambda_0$ is 
nondecreasing along the Ricci flow [P1], we obtain 
\ba
\|u\|^2_{\frac{2n}{n-2}} \le c(1+\frac{1}{\lambda_0(g_0)} (1-\frac{\min_{0} R^-}{4})) \left(\int_M (|\nabla u|^2+\frac{R}{4}u^2)dvol \right)
\ea
which leads to (\ref{sobolevD}). \qed \\

\noindent {\bf Proof of Theorem $\mbox{D}^*$} This is similar to the above proof. \qed \\

\sect{The $\kappa$-noncollapsing estimate}

It is obvious that Theorem E and Theorem $\mbox{E}^*$ follow from Theorem D, 
Theorem $\mbox{D}^*$ and the following lemma. 

\begin{lem} Consider the Riemannian manifold $(M,g)$ for a given metric $g$, such that 
for some $A>0$ and $B>0$ the Sobolev inequality 
\ba \label{lemmasobolev}
\left( \int_M |u|^{\frac{2n}{n-2}} dvol \right)^{\frac{n-2}{n}} \le 
A\int_M (|\nabla u|^2+\frac{R}{4}u^2) dvol +B\int_M u^2 dvol
\ea
holds true for all $u \in W^{1,2}(M)$. 
Let $L>0$. Assume $R\le \frac{1}{r^2}$ on a geodesic ball $B(x, r)$ with $0<r\le L$. Then 
there holds 
\ba \label{noncollapse1}
vol(B(x, r)) \ge \left(\frac{1}{2^{n+3}A+2BL^2}\right)^{\frac{n}{2}} r^n.
\ea
\end{lem} 
\Pf  Let $L>0$.  Assume that 
$R\le \frac{1}{r^2}$ on a closed geodesic ball $B(x_0, r)$ with $0<r\le L$, but the estimate 
(\ref{noncollapse1}) does not hold, i.e. 
\ba \label{small1}
vol(B(x_0, r)) < \delta r^n, 
\ea
where  
\ba
\delta=\left(\frac{1}{2^{n+3}A+2BL^2}\right)^{\frac{n}{2}}. 
\ea
We derive a contradiction. Set 
$\bar g=\frac{1}{r^2} g$. Then we have for $\bar g$
\ba \label{small2}
vol(B(x_0, 1)) < \delta
\ea
and $R\le 1$ on $B(x_0, 1)$.
Moreover, (\ref{lemmasobolev}) leads to the following 
Sobolev inequality for $\bar g$
\ba \label{newsobolev2}
\left( \int_M |u|^{\frac{2n}{n-2}}  \right)^{\frac{n-2}{n}} \le 
A\int_M (|\nabla u|^2+\frac{R}{4}u^2)  +BL^2\int_M u^2,
\ea
where the notation of the volume form is omitted. 
For $u \in C^{\infty}(M)$ with support contained in 
$B(x_0, 1)$ we then have 
\ba \label{ballsobolev1}
\left( \int_{B(x_0, 1)} |u|^{\frac{2n}{n-2}} \right)^{\frac{n-2}{n}} \le 
A\int_{B(x_0, 1)} (|\nabla u|^2+\frac{1}{4}u^2)  +BL^2\int_{B(x_0,1)} u^2.
\ea 
By H\"{o}lder's inequality and (\ref{small2}) we have
\ba
\int_{B(x_0, 1)} u^2  \le \delta^{\frac{2}{n}} \left(\int_{B(x_0, 1)} |u|^{\frac{2n}{n-2}}\right)^{\frac{n-2}{n}}.
\ea
Hence we deduce
\ba \label{ballsobolev2}
\left( \int_{B(x_0, 1)} |u|^{\frac{2n}{n-2}} \right)^{\frac{n-2}{n}} &\le& 
A\int_{B(x_0, 1)} |\nabla u|^2  +(\frac{A}{4}+BL^2)\delta^{\frac{2}{n}} \left(\int_{B(x_0, 1)} |u|^{\frac{2n}{n-2}}\right)^{\frac{n-2}{n}} \nonumber \\
&\le& A\int_{B(x_0, 1)} |\nabla u|^2 +\frac{1}{2}\left(\int_{B(x_0, 1)} |u|^{\frac{2n}{n-2}}\right)^{\frac{n-2}{n}}.
\ea
It follows that 
\ba \label{ballsobolev3}
\left( \int_{B(x_0, 1)} |u|^{\frac{2n}{n-2}} \right)^{\frac{n-2}{n}} &\le& 2A\int_{B(x_0, 1)} |\nabla u|^2.
\ea

Next consider an arbitrary domain $\Omega \subset B(x_0, 1)$. For $u \in C^{\infty}(\Omega)$ with 
support contained in $\Omega$ we deduce from  (\ref{ballsobolev3}) via H\"older's inequality
\ba
\int_{B(x_0, 1)} |u|^2 \le 2Avol(\Omega)^{\frac{2}{n}}\int_{\Omega} |\nabla u|^2.
\ea
Hence we arrive at the following Faber-Krahn inequality
\ba \label{FK}
\lambda_1(\Omega) vol(\Omega)^{\frac{2}{n}} \ge  \frac{1}{2A},
\ea
where $\lambda_1(\Omega)$ denotes the first Dirichlet eigenvalue of 
$-\Delta$ on $\Omega$. By the proof of [C, Proposition 2.4] in [C] we then 
infer
\ba \label{large1}
vol(B(x,\rho)) \ge \left(\frac{1}{2^{n+3}A}\right)^{\frac{n}{2}} \rho^n
\ea
for all $B(x, \rho) \subset B(x_0, 1)$. Consequently we have 
\ba
vol(B(x_0, 1)) \ge \left(\frac{1}{2^{n+3}A}\right)^{\frac{n}{2}},
\ea
contradicting (\ref{small2}).  

For the convenience of the reader, we reproduce here the arguments in the proof 
of [C, Proposition 2.4] in [C].  Consider $B(x, \rho) \subset B(x_0, 1)$. 
Set $u(y)=\rho-d(x, y)$. Then we obtain 
\ba
\lambda_1(B(x, \rho)) \equiv \lambda_1( int \, B(x, \rho))\le \frac{vol(B(x, r))}{\int_{B(x, \rho/2)} u^2} \le 
\frac{4vol(B(x, \rho))}{\rho^2 vol(B(x, \rho/2))}.
\ea
By (\ref{FK}) we then infer 
\ba \label{iteration1}
vol(B(x, \rho)) \ge \left(\frac{\rho^2}{2A}\right)^{\frac{n}{n+2}} 4^{-\frac{n}{n+2}}
vol(B(x, \frac{\rho}{2}))^{\frac{n}{n+2}}.
\ea
Iterating (\ref{iteration1}) we obtain 
\ba
vol(B(x, \rho)) \ge \left(\frac{\rho^2}{2A}\right)^{\sum_{l=1}^m (\frac{n}{n+2})^l} 4^{-\sum_{l=1}^m l(\frac{n}{n+2})^l}
vol(B(x, \frac{\rho}{2^m}))^{(\frac{n}{n+2})^m}
\ea
for all natural numbers $m\ge 1$. Letting $m\rightarrow \infty$ we finally arrive at 
\ba
vol(B(x, \rho)) &\ge& \left(\frac{\rho^2}{2A}\right)^{\sum_{l=1}^{\infty} (\frac{n}{n+2})^l} 4^{-\sum_{l=1}^{\infty} l(\frac{n}{n+2})^l}\nonumber \\
&=&\left(\frac{\rho^2}{2A}\right)^{\frac{n}{2}} 4^{-\frac{n(n+2)}{4}}=\left(\frac{1}{2^{n+3}A}\right)^{\frac{n}{2}} \rho^n.
\ea
\qed \\

\noindent {\bf Proof of Theorem F} This theorem follows from Theorem E and Cheeger-Gromov-Hamilton compactness 
theorem.   \qed

\vspace{2cm}

\noindent {\bf {\Large Appendices}} \\

\appendix
\sect{The logarithmic Sobolev inequalities on the Euclidean space}

In this appendix we review several versions of the logarithmic Sobolev inequality on the 
euclidean space for the purpose of presenting the background of  the logarithmic 
Sobolev inequalitites.  These versions are equivalent to each other.\\

\noindent {\bf 1. The Gaussian version} \\

This is the original version of L.~Gross.

\begin{theo}  Let $u\in W^{1,2}_{loc}({\bf R}^n)$ satisfy 
 $\int_{{\bf R}^n} u^2d\mu=1$, where
 \be
d\mu=(2\pi)^{-\frac{n}{2}}e^{-\frac{|x|^2}{2}}dx.
\tag{A.1}
\ee
 Then 
\be
\int_{{\bf R}^n} u^2 \ln u^2 d\mu \le 2 \int_{{\bf R}^n} |\nabla u|^2 d\mu.
\tag{A.2}
\ee
\end{theo}

\vspace{1cm}

\noindent {\bf 2. The straight (Euclidean volume element)   version} \\

\begin{theo}  There holds  
\be
\int u^2 \ln u^2 dx \le 2\int |\nabla u|^2 dx,
\tag{A.3}
\ee
provided that $u \in W^{1,2}({\bf R}^n)$  and $\int u^2 dx=(2 \pi)^{n/2} e^n.$  Equivalently, for $\beta>0$, 
\be
\int u^2 \ln u^2 dx \le 2\int |\nabla u|^2dx + \beta \ln \beta-\frac{n}{2} \beta\ln (2 \pi e^2), 
\tag{A.4}
\ee
provided that $u \in W^{1,2}({\bf R}^n)$  and $\int u^2=\beta$.
\end{theo}

\vspace{1cm}

\noindent {\bf 3. The log gradient version} \\

It appears to be stronger than the other versions because of the logarithm in front of the 
Dirichlet integral of $u$. 

\begin{theo}  There holds
\be
\int u^2 \ln u^2 dx \le \frac{n}{2} \ln \left[ \frac{2}{\pi n e } \int |\nabla u|^2 dx \right],
\tag{A.5}
\ee
provided that $u \in W^{1,2}({\bf R}^n)$  and $\int u^2 dx =1$. \\
\end{theo}
 
\vspace{1cm}
 
\noindent {\bf 4. The entropy version (as formulated in [P1])} \\

This version is intimately related to Perelman's entropy functional ${\mathcal W}$. Indeed, it can be 
viewed as the motivation for $\mathcal W$. 

\begin{theo}  There holds 
\be
\int (\frac{1}{2}|\nabla f|^2+f-n)e^{-f}dx \ge 0,
\tag{A.6}
\ee
provided that $f \in W^{1,2}_{loc}({\bf R}^n)$ and $\int e^{-f}dx=(2\pi)^{n/2}$. 
\end{theo}

\sect{The estimate for $e^{tH}$}

In this appendix we present the proof of Theorem \ref{D1}.  The global case $\sigma^*=\infty$ of this theorem follows from  
Corollary 2.2.8 in [D].  On the other hand,  
the proof of this corollary in [D] can easily be extended to cover  the local case $\sigma^*<\infty$, as is done 
below.    The global case is customarily phrased in terms of ``ultracontractivity", i.e. the logarithmic Sobolev inequality implies the ultracontractivity of $e^{-tH}$, see e.g. [D].  Note that the global 
case suffices for the main purpose of this paper. The local case should be useful for further applications.
\\

\noindent {\bf Proof of Theorem \ref{D1}} \\

{\bf Part 1}  We first assume $\Psi \ge 0$, i.e. $\min \Psi^-=0$. It follows from (\ref{Dlog1})
\be \label{Dlog2}
\int_M u^2 \ln u^2  \le \sigma Q(u)+ \beta(\sigma)\|u\|_2^2+\|u\|_2^2 \ln \|u\|_2^2
\tag{B.1}
\ee
for all $u \in W^{1,2}(M)$. Here the notation of the volume form is omitted.  Replacing $u$ by $|u|^{p/2}$ for $p>2$ and 
$u \in W^{1,2}(M) \cap L^{\infty}(M)$ we deduce 
\be \label{Dlog3}
p \int_M |u|^p \ln |u|  \le \sigma Q(|u|^{\frac{p}{2}})+ \beta(\sigma) \|u\|_p^p
+p \|u\|_p^p \ln \|u\|^2_p.
\tag{B.2}
\ee
Since 
\be
Q(|u|^{\frac{p}{2}})=\frac{p^2}{4(p-1)} Q(|u|, |u|^{p-1})
\tag{B.3}
\ee
we arrive at 
\be \label{Dlog4}
\int_M |u|^p \ln |u| \le \frac{\sigma p}{4(p-1)} Q(|u|, |u|^{p-1})+ \frac{\beta(\sigma)}{p} \|u\|_p^p
+\|u\|_p^p \ln \|u\|_p. 
\tag{B.4}
\ee
By the nonincreasing property of $\beta$ we then infer, replacing $\sigma$ by 
$\frac{4(p-1)}{p} \sigma$
\be \label{Dlog5}
\int_M |u|^p \ln |u| \le  \sigma Q(|u|, |u|^{p-1})+ \frac{\beta(\sigma)}{p} \|u\|_p^p
+\|u\|_p^p \ln \|u\|_p 
\tag{B.5}
\ee
for $\sigma \in (0, \frac{p}{4(p-1)}\sigma^*]$. 
\\

\noindent {\bf Part 2} We continue with the assumption $\Psi \ge 0$. Consider $0<t \le \frac{1}{4} \sigma^*$.  
Let $\sigma(p)$ be a nonnegative continuous function for $p\ge 2$ such that $\sigma(p) \in (0,  
\frac{p}{4(p-1)}\sigma^*]$ for $p>2$, which will be chosen later. Then we have 
\be \label{Dlog6}
\int_M |u|^p \ln |u| \le  \sigma(p) Q(|u|, |u|^{p-1})+ \Gamma(p) \|u\|_p^p
+\|u\|_p^p \ln \|u\|_p 
\tag{B.6}
\ee
for each $p>2$ and all $u \in W^{1,2}(M) \cap L^{\infty}(M)$, 
where $\Gamma(p)=\frac{\beta(\sigma(p))}{p}.$  Define  the function $p(s)$ 
for $0\le s <t$ by 
\be
\frac{dp}{ds}=\frac{p}{\sigma(p)}, p(0)=2.
\tag{B.7}
\ee
Assume that 
\be \label{infty}
p(s) \rightarrow \infty
\tag{B.8}
\ee
as $s \rightarrow t$. We also define the function $N(s)$ for $0\le s<t$ by 
\be
\frac{dN}{ds}=\frac{\Gamma(p(s))}{\sigma(s)}, N(0)=0
\tag{B.9}
\ee
and set 
\be
N^*=\lim\limits_{s\rightarrow t} \equiv \int_2^{\infty} \frac{\Gamma(p)}{p}dp.
\tag{B.10}
\ee

For $u\in W^{1,2}(M)\cap L^{\infty}(M)$ with $u\ge 0$ we set $u_s=e^{-sH}u$ for $0<s<t$.  
By the contraction properties of $e^{-sH}$ we have $u_s \in W^{1,2}(M) \cap L^{\infty}(M)$ for 
all $s$. If $\Psi \in C^{\infty}(M)$ we have 
for a fixed $q>2$ 
\be
\frac{d}{ds}\|u_s\|_q^q =q \int_M \frac{\partial u_s}{\partial s} \cdot u_s^{q-1} =
-q \int_M Hu_s \cdot u_s^{q-1}.
\tag{B.11}
\ee
Hence
\be
\frac{d}{ds}\|u_s\|_q^q =-qQ(u_s, u_s^{q-1}).
\tag{B.12}
\ee
In the general case $\Psi\in L^{\infty}(M)$, this formula follows from the spectral formula for 
$e^{-sH}$.  
Using this formula we compute  
\be
\frac{d}{ds}\ln (e^{-N(s)}\|u_s\|_{p(s)})=\frac{d}{ds}\left(-N(s)+\frac{1}{p(s)} 
\ln \|u_s\|_{p(s)}^{p(s)}\right) \nonumber
\ee
\be
=\frac{\Gamma}{\sigma}-\frac{1}{p^2}\frac{p}{\sigma} \ln \|u_s\|_p^p +\frac{1}{p} \|u_s\|_p^{-p}
\left(-pQ(u_s, u_s^{p-1})+\frac{p}{\sigma} \int_M u_s^p \ln u_s\right) \nonumber 
\ee
\be
=\frac{1}{\sigma}\|u_s\|^{-p}_p \left( \int_M u_s^p \ln u_s -
\sigma Q(u_s, u_s^{p-1})-\Gamma \|u_s\|_p^p-\|u_s\|_p^p \ln \|u_s\|_p \right).
\tag{B.13}
\ee
By (\ref{Dlog6}) this is nonpositive. Hence $e^{-N(s)}\|u_s\|_{p(s)}$ is 
nonincreasing, which leads to 
\be
\|e^{-sH}u\|_{p(s)} \le e^{N(s)} \|f\|_2
\tag{B.14}
\ee
for
 all $0\le s<t$.
By the contraction properties we have $\|e^{-tH}u\|_{p(s)}  \le \|e^{-sH}u\|_{p(s)}$, whence
\be
\|e^{-tH}u\|_{p(s)}  \le  e^{N(s)} \|f\|_2
\tag{B.15}
\ee
for all $0\le s<t$. It follows that
\be
\|e^{-tH}u\|_{\infty} \le e^{N^*} \|u\|_2.
\tag{B.16}
\ee
This estimate extends to $u\in L^2(M)$ with $u\ge 0$ by an approximation. For a 
general $u \in L^2(M)$ we use the pointwise inequality $|e^{-tH}u| \le e^{-tH}|f|$ 
(a consequence of the positivity preserving property) to deduce
\be
\|e^{-tH}u\|_{\infty}\le \|e^{-tH}|u|\|_{\infty} \le e^{N^*} \|u\|_2.
\tag{B.17}
\ee

Now we choose  
\be \tag{B.18}
\sigma(p)=\frac{2t}{p}
\ee
for $p \ge 2$. Then $p(s)=\frac{2t}{t-s}$.  One readily sees that $\sigma(p)\in (0, \frac{p}{4(p-1)} \sigma^*]$ 
for $p>2$ and $p(s)\rightarrow \infty$ as $s \rightarrow t$. We have for this choice
\be
N^*=\frac{1}{2t} \int_0^t \beta(\sigma)d\sigma.
\tag{B.19}
\ee
Hence we arrive at 
\be \label{positivePsi}
\|e^{-tH}u\|_{\infty} \le e^{\tau(t)} \|u\|_2
\tag{B.20}
\ee
for all $u \in L^2(M)$ and  $0<t\le \frac{1}{4}\sigma^*$. 
\\

\noindent {\bf Part 3} For a general $\Psi$, we consider  $\bar \Psi=\Psi-\min \Psi^-$ 
and denote the corresponding operator and quadratic form by $\bar H$ and 
$\bar Q$ respectively. We have by (\ref{Dlog1})
\be \label{Dlog7}
\int_M u^2 \ln u^2 dvol \le \sigma \bar Q(u)+ \bar \beta(\sigma)
\tag{B.21}
\ee
for all $u \in L^2(M)$ with $\|u\|_2=1$,
where $\bar \beta(\sigma)=\beta(\sigma)+\sigma \min \Psi^-.$ 
We apply  (\ref{positivePsi}) to deduce for $0<t\le \frac{1}{4}\sigma^*$ and 
$u \in L^2(M)$
\be 
\|e^{-t\bar H} u\|_{\infty} \le e^{\frac{1}{2t} \int_0^t \bar \beta(\sigma) d\sigma}  \|u\|_2 
= e^{\tau(t)+\frac{t}{4}  \min \Psi^-} \|u\|_2. \tag{B.22}
\ee
The desired estimate (\ref{heat1}) follows. 

The estimate (\ref{heat2}) follows from 
(\ref{heat1}) in terms of duality. Namely we have for 
$u, v \in L^2(M)$
\be
\int_M v e^{-tH}u = \int_M u e^{-tH}v 
\le \|e^{-tH}v\|_{\infty} \|u\|_1 \le  e^{\tau(t)-\frac{3t}{4}  \min \Psi^-}\|v\|_2 \|u\|_1
\tag{B.23}
\ee 
It follows that 
\be
\|e^{-tH}u\|_2 \le e^{\tau(t)-\frac{3t}{4}  \min \Psi^-}\|u\|_1
\tag{B.24}
\ee
and then 
\be
\|e^{-tH}u\|_{\infty} \le e^{\tau(\frac{t}{2})-\frac{3t}{8}  \min \Psi^-}\|e^{-\frac{t}{2}H}u\|_2 \le
e^{2\tau(\frac{t}{2})-\frac{3t}{4}  \min \Psi^-} \|u\|_1.
\tag{B.25}
\ee
By Lemma \ref{extensionlemma}, we arrive at (\ref{heat2}) for all $u \in L^1(M)$.
The estimate (\ref{heat2}) also follows from the arguments in Part 2 by choosing 
$\sigma(p)=\frac{t}{p}$ and $p(s)=\frac{t}{t-s}$.  

\sect{From the estimate for $e^{-tH}$ to the Sobolev inequality}

In this appendix we present the proof of Theorem \ref{D2}. We also present a more general 
result Theorem \ref{general}, and its implication for the Ricci flow.  
Consider a compact Riemannian manifold $(M, g)$ of dimension $n \ge 1$ and 
$\Psi \in L^{\infty}(M)$ as in the set-up 
for Theorem \ref{D2}.   If $Q \ge 0$, then we define the spectral square root
$H^{\frac{1}{2}}$ of the operator $H=-\Delta+\Psi$ as follows. For 
$u=\sum_{i \ge 1} a_i \phi_i \in L^2(M)$ we set 
\ba
H^{\frac{1}{2}}u=\sum_{i\ge 1} \lambda_i^{\frac{1}{2}} a_i \phi_i,
\ea
whenever the series converges in $L^2(M)$. 

\begin{lem} \label{Qlemma} Assume $Q \ge 0$. Then $H^{\frac{1}{2}}$ is a bounded operator from  $W^{1,2}(M)$ to $L^2(M)$.
Indeed there holds for all $u \in W^{1,2}(M)$
\ba \label{HQ}
\|H^{\frac{1}{2}}u\|_2^2=Q(u).
\ea
\end{lem}
\Pf  For $u=\sum_{\i \ge 1} a_i \phi_i \in C^2(M)$ there holds 
$Q(u)=<Hu, u>_2 = \sum_{i\ge 1} \lambda_i a_i^2.$
By approximation, we derive  $Q(u)=\sum_{i \ge 1} \lambda_i a_i^2$ for all 
$u \in W^{1,2}(M)$.  Now we have for $N \ge 1$ 
\ba
\|\sum_{1\le i \le N} \lambda_i^{\frac{1}{2}} a_i \phi_i \|^2_2
=\sum_{1\le i \le N} \lambda_i a_i^2.
\ea
Taking the limit as $N \rightarrow \infty$ we infer 
$\|H^{\frac{1}{2}}u\|^2_2 = Q(u).$
\qed \\

If $Q>0$, i.e. the first eigenvalue of $H$ is positive, then the inverse $H^{-\frac{1}{2}}:
L^2(M) \rightarrow W^{1,2}(M)$ of $H^{\frac{1}{2}}$ exists. We have
$H^{-\frac{1}{2}}u=\sum_{i\ge 1} \lambda_i^{-\frac{1}{2}} a_i \phi_i$
for $u=\sum_{i\ge 1} a_i \phi_i \in L^2(M)$.
More generally, we define $H^{-\frac{1}{2}}$ in the case $Q \ge 0$ by 
$H^{-\frac{1}{2}}u=\sum_{\lambda_i>0} \lambda_i^{-\frac{1}{2}} a_i \phi_i$ for 
$u=\sum_{i\ge 1} a_i \phi_i \in L^2(M)$. 

\begin{lem} \label{halflemma} Assume $Q \ge 0$. We set $\phi_1^*=\phi_1$ if 
$\lambda_1=0$ and $\phi_1^*=0$ if $\lambda_1>0$.  There holds 
\ba \label{halfformula}
H^{-\frac{1}{2}} u =\Gamma(\frac{1}{2})^{-1} \int_0^{\infty} t^{-\frac{1}{2}} 
e^{-tH} u dt
\ea
for all $u \in L^2(M)$ with $u \perp \phi_1^*$.  Moreover, if $u \in L^2(M)$ with $u \perp 
\phi_1^*$ satisfies  $\|e^{-tH}u\|_{\infty} \le \phi(t)$ on an open interval 
$(a, b) \subset (0, \infty)$ for a 
nonnegative continuous function $\phi$, then there holds 
\ba \label{infinity}
\|\int_a^b t^{-\frac{1}{2}} 
e^{-tH} u dt\|_{\infty} \le \int_a^b t^{-\frac{1}{2}} 
\phi(t) dt.
\ea
\end{lem} 
\Pf  For $u\in L^2(M)$ with $u \perp \phi_1^*$ we write $u =\sum_{\lambda_i>0} a_i \phi_i$, where the series converges in 
$L^2(M)$.  We have $e^{-tH}u=\sum_{\lambda_i>0} e^{-\lambda_i t} a_i \phi_i$.
We have 
\ba
 \sum_{\lambda_i>0} \int_0^{\infty} t^{-\frac{1}{2}} e^{-\lambda_i t} a_i \phi_i dt
=\Gamma(\frac{1}{2}) \sum_{\lambda_i>0} \lambda_i^{-\frac{1}{2}} a_i \phi_i=\Gamma(\frac{1}{2})H^{-\frac{1}{2}}u.
\ea
Hence the formula (\ref{halfformula}) follows.
Next we note that convergence in $L^2$ implies 
almost everywhere convergence. Moreover, if $u_k$ converges to $u$ almost everywhere,
then $\|u\|_{\infty} \le \liminf \|u_k\|_{\infty}$.   These two facts lead to 
(\ref{infinity}).  \qed \\

Next we recall, for the sake of clarity and precise estimates, the Marcinkiewicz interpolation theorem [Sa] and the Riesz-Thorin 
interpolation theorem [Sa], which we formulate in the 
special case of the measure space $(M, \mu)$, where $\mu$ denotes the Lebesgue measure 
associated with the volume element $dvol$ of $g$.    

\begin{theo} \label{Mar} (Marcinkiewicz interploation theorem) \, \, Let $L$ be an additive operator from $L^{\infty}(M)$ to the 
space of measurable functions on $M$.  Let $1 \le p_0 \le q_0 \le \infty$ and  
$1 \le p_1 \le q_1 \le \infty$ with $q_0 \not =q_1$. Assume that $L$ is of weak type 
$(p_0, q_0)$ with constant $K_0$ and of weak type $(p_1, q_1)$ with 
constant $K_1$, i.e. 
\ba
\mu(\{|L(u)| >\alpha\}) \le \left(K_0 \frac{\|u\|_{p_0}}{\alpha}\right)^{q_0}
\ea
and 
\ba
\mu(\{|L(u)| >\alpha\}) \le \left(K_1 \frac{\|u\|_{p_1}}{\alpha}\right)^{q_1}
\ea
for all $u \in L^{\infty}(M)$.
Then $L$ is of  type $(p_t, q_t)$ on $L^{\infty}(M)$ with constant $K_t$ for each $0<t<1$, i.e. 
\ba
\|L(u)\|_{q_t} \le K_t \|u\|_{p_t}
\ea
for all $u \in L^{\infty}(M)$ and $\alpha>0$, where   
\ba \label{Lpq}
\frac{1}{p_t}=\frac{1-t}{p_0}+\frac{t}{p_1}, \frac{1}{q_t}=\frac{1-t}{q_0}+
\frac{t}{q_1},
\ea
\ba \label{Kbound1}
K_t\le KK_0^{1-t}K_1^t,
\ea
and $K=K(p_0, q_0, p_1, q_1, t)$ is bounded for $0<\epsilon\le t \le 1-\epsilon$ with 
each given $\epsilon>0$, but tends to infinity as $t\rightarrow 0$ or $t \rightarrow 1$. 

It follows that for each $0<t<1$, $L$ extends uniquely to an additive operator $L: L^{p_t}(M) 
\rightarrow L^{q_t}(M)$ with the bound (\ref{Lpq}).
\end{theo} 

This follows from [Sa, Theorem 5.2]. The space of simple functions is used 
in [Sa, Theorem 5.2] instead of $L^{\infty}(M)$. Moreover, $L$ is only assumed to be sublinear. 
Note that Theorem \ref{Mar} holds both in the set-up of real-valued functions and the set-up of 
complex-valued functions.

\begin{theo} \label{RT} (Riesz-Thorin interpolation theorem) \, \, Let $L$  be a linear operator from 
$L^{\infty}_{{C}}(M)$, i.e. the complex-valued $L^{\infty}(M)$,  to the space of 
complex valued measurable functions on $M$. Let $1\le p_0, p_1, q_0, q_1 \le \infty$. 
Assume that $L$ is of type $(p_0, q_0)$ on $L^{\infty}_C(M)$ with constant $K_0$, and of type $(p_1, q_1)$ 
on $L^{\infty}_C(M)$ with constant 
$K_1$. Then $L$ is of type $(p_t, q_t)$ on $L^{\infty}_C(M)$ with constant $K_t$ for each 
$0\le t \le 1$, where $p_t$ and $q_t$ are given by (\ref{Lpq}) and 
\ba \label{Kbound2}
K_t \le K_0^{1-t} K_1^t.
\ea
Consequently, for each $0\le t \le 1$, $L$ extends uniquely to a linear operator 
$L: L^{p_t}_C(M) \rightarrow L^{q_t}_C(M)$ with the bound 
\ba
\|L(u)\|_{q_t} \le K_t \|u\|_{p_t}
\ea
for all $u \in L^{p_t}_C(M)$, where $L^p_C(M)$ denotes the complex-valued 
$L^p(M)$.  

If we replace the complex-valued spaces by real-valued spaces, then the same holds except that 
(\ref{Kbound2}) is replaced by 
\ba \label{Kbound3}
K_t \le 2 K_0^{1-t} K_1^t.
\ea
The bound (\ref{Kbound2}) still holds in the set-up of real-valued functions, provided that 
$p_0 \le q_0, p_1 \le q_1$, or $T$ is a positive operator.
\end{theo}

Now we are ready to prove Theorem \ref{D2}. We present a more general result which implies 
Theorem \ref{D2}.

\begin{theo}  \label{general} 1) Let $\mu>1$. Assume that $\Psi\ge 0$ and  
for some $c>0$ the inequality 
\ba \label{heatcondition}
\|e^{-tH}u\|_{\infty} \le c t^{-\frac{\mu}{4}} \|u\|_2
\ea 
holds true for each $t>0$ and all $u\in L^2(M)$.  Let $1<p<\mu$.  Then  
there holds 
\ba \label{Hestimate1}
\|H^{-\frac{1}{2}}u\|_{\frac{\mu p}{\mu-p}} \le C(c, \mu, p) \|u\|_p  
\ea
for all $u \in L^p(M)$, where the positive constant 
 $C(\mu, c, p)$ can be bounded from above in terms of upper bounds for $c$, $\mu$,  
$\frac{1}{\mu-p}$ and $\frac{1}{p-1}$.  Consequently, there holds
\ba \label{Hestimate1a}
\|u\|_{\frac{\mu p}{\mu-p}} \le C(c, \mu, p) \|H^{\frac{1}{2}}u\|_p
\ea
for all $u \in W^{1,p}(M)$.
\\
2) Let $\mu>1$. Assume that for some $c>0$ the inequality  
\ba 
\|e^{-tH}u\|_{\infty} \le c t^{-\frac{\mu}{4}}\|u\|_2
\ea 
holds true for each $0<t<1$ and all $u\in L^2(M)$. Set $H_0=H-\inf \Psi^-+1$.  Let $1<p<\mu$.  Then there  
holds
\ba \label{Hestimate2}
\|H_0^{-\frac{1}{2}}u\|_{\frac{\mu p}{\mu-p}} \le C(\mu,c, p) \|u\|_p
\ea
for all $u \in L^p(M)$, where the positive constant $C(\mu,c, p)$ has the same property as 
the $C(\mu, c, p)$ above. Consequently, there holds 
\ba \label{Hestimate2a}
\|u\|_{\frac{\mu p}{\mu -p}} \le C(\mu, c, p) \|H_0^{\frac{1}{2}}u\|_p
\ea
for all $u \in W^{1,p}(M)$.
\end{theo} 
\Pf 1) For simplicity, we work in the set-up of real-valued functions.   The case $p=2$ follows from  [D, Theorem 2.4.2]. 
The proof of that theorem in [D] 
extends in a standard way to the general case of (\ref{Hestimate1}), so we follow it here.  By the proof of Theorem \ref{D1} in Appendix B we have 
with $\tau(t)=\ln c_1-\frac{\mu}{4} \ln t$
\ba
\|e^{-tH}u\|_{\infty} \le e^{2 \tau(\frac{t}{2})}\|u\|_1 = \frac{2^{\frac{\mu}{2}}c^2}{t^{\frac{\mu}{2}}}
\|u\|_1
\ea
for all $u \in L^1(M)$ and all $t>0$. 
On the other hand, we have by Lemma \ref{contractionlemma} $\|e^{-tH}u\|_{\infty} \le \|u\|_{\infty}$
for all $u \in L^{\infty}(M)$.  By Theorem \ref{RT} we then have
\ba \label{inf-p}
\|e^{-tH}u\|_{\infty} \le  \left(\frac{2^{\frac{\mu}{2}}c^2}{t^{\frac{\mu}{2}}}\right)^{\frac{1}{p}} 
\|u\|_p
\ea
for each $1 \le p \le \infty$ and all $u \in L^p(M)$. 

Next we consider $1 \le p<\mu$ and set
\ba
\frac{1}{q}=\frac{1}{p}-\frac{1}{\mu}, \,\,\, \mbox{i.e.} \, \, \, q=\frac{\mu p}{\mu-p}.
\ea 
Observe that (\ref{heatcondition}) implies that the first eigenvalue $\lambda_1$ of $H$ is positive. 
Otherwise, since $Q \ge 0$, $\lambda_1$ would be zero. Then 
$e^{-tH}\phi_1=\phi_1$ for all $t>0$.  This contradicts 
(\ref{heatcondition}).  Thus (\ref{halfformula}) is valid for all $u \in L^2(M)$.  We show that $H^{-\frac{1}{2}}:L^{\infty}(M) \rightarrow 
L^2(M)$ is of weak type $(p, q)$. For a given $T \in (0, \infty)$ we write 
$H^{-\frac{1}{2}}(u)=G_{0, T}(u)+G_{T,\infty}(u)$, where 
\ba
G_{a, b}(u)=\Gamma(\frac{1}{2})^{-1} \int_a^b t^{-\frac{1}{2}} e^{-tH}udt.
\ea
We have by Lemma \ref{halflemma} and (\ref{inf-p}) 
\ba
\|G_{T, \infty}(u)\|_{\infty} &\le& \Gamma(\frac{1}{2})^{-1}\int_T^{\infty}
t^{-\frac{1}{2}} \left(\frac{2^{\frac{\mu}{2}}c_1^2}{t^{\frac{\mu}{2}}}\right)^{\frac{1}{p}}  \|u\|_p dt 
\nonumber \\ &=& c_1(\mu, p) T^{\frac{1}{2}-\frac{\mu}{2p}}\|u\|_p 
\ea
for all $u \in L^{p}(M)$, where 
\ba
c_1(\mu, p)=\frac{2p}{\mu-p}\Gamma(\frac{1}{2})^{-1} \left(2^{\frac{\mu}{2}}c^2\right)^{\frac{1}{p}}.
\ea
On the other hand, we have by Lemma \ref{contractionlemma}
\ba
\|G_{0, T}(u)\|_{p} \le  \Gamma(\frac{1}{2})^{-1} \int_0^T t^{-\frac{1}{2}} \|u\|_pdt 
=2 \Gamma(\frac{1}{2})^{-1} T^{\frac{1}{2}} \|u\|_p 
\ea
for all $u \in L^p(M)$.  Given $u \in L^{\infty}(M)$ and $\alpha>0$ we define $T$ by 
\ba
\frac{\alpha}{2}=c_1(\mu, p) \|u\|_p T^{\frac{1}{2}-\frac{\mu}{p}}.
\ea
Then $\|G_{T, \infty}(u)\|_{\infty} \le \frac{\alpha}{2}$ and hence 
\ba \label{weak-p-q}
\mu(\{|H^{-\frac{1}{2}}(u)|>\alpha\}) &\le& \mu(\{|G_{0,T}(u)|>\frac{\alpha}{2}\}) 
\le \left(\frac{\alpha}{2}\right)^{-p} \|G_{0, T}(u)\|_p^p \nonumber \\
&\le& \left(\frac{\alpha}{2}\right)^{-p}\left(2 \Gamma(\frac{1}{2})^{-1} T^{\frac{1}{2}}\right)^p \|u\|_p^p
\nonumber \\
&=& c_2(\mu, p)\left( \frac{\|u\|_p}{\alpha}\right)^q, 
\ea
where 
\ba
c_2(\mu, p)=2^{\frac{p(2\mu-p)}{\mu-p}}c_1(\mu, p)^{\frac{p^2}{\mu-p}} \Gamma(\frac{1}{2})^{-p}.
\ea
It follows that $H^{-\frac{1}{2}}$ is of weak type $(p, q)$ with constant $c_2(\mu, p)^{1/q}$.

Given $1<p<\mu$, we set $\gamma=\max\{\frac{p}{p-1}, 2, \frac{2\mu-p}{\mu-p}\}+1$, $p_0=\frac{\gamma-1}{\gamma}p$ and 
$p_1=\frac{\gamma-1}{\gamma-2}p$. Then $1<p_0<p_1<\mu$ and 
\ba
\frac{1}{p_0}+\frac{1}{p_1}=\frac{2}{p},\, \, \, 
\frac{1}{q_0}+\frac{1}{q_1}=\frac{2}{q},
\ea
where ${1}/{q_0}={1}/{p_0}-{1}/{\mu}$ and ${1}/{q_1}={1}/{p_1}-{1}/{\mu}$, and $q$ is the same as before, i.e. 
$1/q=1/p-1/\mu$.
Applying (\ref{weak-p-q}) and Theorem \ref{Mar} with $t=\frac{1}{2}$ we then arrive at 
(\ref{Hestimate1}) with
\ba
C(\mu, c, p)=K(p_0, q_0, p_1, q_1, \frac{1}{2}) c_2(\mu, p_0)^{\frac{1}{2q_0}} c_2(\mu, p_1)^{\frac{1}{2q_1}}.
\ea
The property of $C(\mu, c, p)$ is easy to see from this formula.

By [Se], the operator $H^{\frac{1}{2}}$ is a pseudo-differential operator of 
 order $1$.  Since $M$ is compact, it follows that $H^{\frac{1}{2}}$ is a bounded operator from $W^{1,p}(M)$ into $L^p(M)$ for all $1<p<\infty$. (The special case 
 $p=2$ is contained in Lemma \ref{Qlemma}).  For $2 \le p <\mu$ (assuming $\mu>2$) we have $W^{1,p}(M) 
 \subset W^{1,2}(M)$, and hence $H^{-\frac{1}{2}}H^{\frac{1}{2}}u=u$ by Lemma 
 \ref{Qlemma}.  
 Replacing $u$ in (\ref{Hestimate1}) by $H^{\frac{1}{2}}u$ for $u \in W^{1, p}(M)$ we 
 then arrive at  (\ref{Hestimate1a}). For $1<p<\min\{2, \mu\}$, we can argue this way to arrive at 
 (\ref{Hestimate1}) for $u \in C^{\infty}(M)$.  By the boundedness of $H^{\frac{1}{2}}: 
 W^{1,p}(M) \rightarrow L^p(M)$ we then arrive at (\ref{Hestimate1}) for all 
 $u \in W^{1,p}(M)$ via approximation. 
 
 By [Se], the operator $H^{-\frac{1}{2}}$ is a pseudo-differential operator of order 
 $-1$. It follows that 
 $H^{-\frac{1}{2}}$ is a bounded map from $W^{1,p}(M)$ into $L^p(M)$ for all 
 $1<p<\infty$.  It also follows that  $H^{-\frac{1}{2}}: L^p(M) \rightarrow 
 W^{1,p}(M)$ is the inverse of $H^{\frac{1}{2}}: W^{1,p}(M) \rightarrow L^p(M)$. 
 Moreover, by approximation the inequality (\ref{Hestimate1a}) also implies the inequality (\ref{Hestimate1}).  \\
2) For $0<t<1$ we have 
for $u\in L^2(M)$
\ba
\|e^{-tH_0}u\|_{\infty} =e^{-t(1-\inf \Psi^-) }\|e^{-tH}u\|_{\infty} \le c t^{-\frac{\mu}{4}} \|u\|_2.
\ea
For $t\ge 1$ we write $t=\frac{m}{2}+t_0$ for a natural number $m$ such that $1/2 \le t_0 <1$. Then 
we have for $u \in L^2(M)$ 
\ba
\|e^{-tH_0}u\|_{\infty}&=&e^{-t(1-\inf \Psi^-)} \|e^{-\frac{m}{2}tH} e^{-t_0H}u\|_{\infty} 
\le e^{-t}\|e^{-t_0H}u\|_{\infty} \nonumber \\
 &\le& c e^{-t} t_0^{-\frac{\mu}{4}} \|u\|_2 \le c2^{\frac{\mu}{4}} e^{-t}\|u\|_2
 \nonumber \\ &\le&
c2^{\frac{\mu}{4}}e^{-\frac{\mu}{4}} \left(\frac{\mu}{4}\right)^{\frac{\mu}{4}} t^{-\frac{\mu}{4}} 
\|u\|_2.
\ea
Hence we can apply the result in 1) to arrive at the desired inequalities 
(\ref{Hestimate2}) and (\ref{Hestimate2a}).  (Note that by the above arguments they  are equivalent to each other.) 
\qed \\

\noindent {\bf Proof of Theorem \ref{D2}} 1) Let $u \in W^{1,2}(M)$.  Applying (\ref{Hestimate1a}) with $p=2$  we arrive 
at 
\ba
\|u\|_{\frac{2\mu}{\mu-2}} \le C(\mu, c, 2)\|H^{\frac{1}{2}}u\|_2.
\ea
Combining this with (\ref{HQ}) we then obtain the desired inequality.  \\
2) This is similar to 1). Note that the quadratic form of $H-\inf \Psi^-+1$ is 
$Q(u)+(1-\inf \Psi^-) \|u\|_2^2$. \qed \\

Combining Theorem A, Theorem C, Theorem \ref{D1} and Theorem \ref{general} 
we obtain the following two results for the Ricci flow, which extend  
Theorem D and Theorem $\mbox{D}^*$.  Let $g=g(t)$ be a smooth soluton of the Ricci flow 
on $M \times [0, T)$ as before.

\begin{theo}  Assume that $R_{g_0} \ge 0$ and $\lambda_0(g_0)>0$.  Let $1<p<n$.   There is a  positive constant $C$ depending only on 
the dimension $n$, a positive lower bound for $\lambda_0(g_0)$, a positive lower bound for $vol_{g_0}(M)$, an upper bound for $C_S(M,g_0)$, an upper bound for 
$\frac{1}{p-1}$, and an upper bound for $\frac{1}{n-p}$, such that for each $t \in [0, T)$ and 
all $u \in W^{1,p}(M)$ there holds 
\ba
\|u\|_{\frac{np}{n-p}} \le C \|(-\Delta+\frac{R}{4})^{\frac{1}{2}}u\|_p.
\ea
\end{theo}

\begin{theo} Assume $T<\infty$ and $1<p<n$.  There is a positive constant $C$ depending only on the dimension $n$, 
a nonpositive lower bound for $R_{g_0}$, a positive lower bound for $vol_{g_0}(M)$, an upper bound for $C_S(M,g_0)$,
an upper bound for $T$, an upper bound for $\frac{1}{p-1}$, and an upper bound 
for $\frac{1}{n-p}$,  such that for each $t \in [0, T)$ and 
all $u \in W^{1,p}(M)$ there holds 
\ba
\|u\|_{\frac{np}{n-p}} \le C \|H_0^{\frac{1}{2}}u\|_p,
\ea
where 
\ba
H_0=-\Delta+\frac{R}{4}-\frac{\min R_{g_0}^-}{4}+1.
\ea
\end{theo}

\end{document}